\documentclass[final]{siamltex} 
\usepackage{amssymb,amsmath,amsfonts,bm,epsfig,epstopdf,eqnarray}
\usepackage{graphicx,psfrag,subfig}
\usepackage{algorithm}
\usepackage{algorithmic}
\usepackage[usenames,dvipsnames]{color}
\usepackage{stmaryrd}
\usepackage{ulem}
\usepackage{marginnote}
\usepackage[displaymath, mathlines]{lineno}
\usepackage{environ}
\usepackage[yyyymmdd,hhmmss]{datetime}

\def\Ab{{\boldsymbol A}}
\def\Bb{{\boldsymbol B}}
\def\Cb{{\boldsymbol C}}
\def\Db{{\boldsymbol D}}

\def\Gb{{\boldsymbol G}}
\def\Hb{{\boldsymbol H}}
\def\Ib{{\boldsymbol I}}

\def\Kb{{\boldsymbol K}}
\def\Lb{{\boldsymbol L}}
\def\Mb{{\boldsymbol M}}

\def\Pb{{\boldsymbol P}}
\def\Qb{{\boldsymbol Q}}
\def\Rb{{\boldsymbol R}}
\def\Sb{{\boldsymbol S}}
\def\Tb{{\boldsymbol T}}
\def\Ub{{\boldsymbol U}}
\def\Vb{{\boldsymbol V}}
\def\Wb{{\boldsymbol W}}
\def\Xb{{\boldsymbol X}}
\def\Yb{{\boldsymbol Y}}

\def\mub{{\boldsymbol         \mu}}

\def\omegab{{\boldsymbol      \omega}}

\def\qb{{\boldsymbol q}}

\def\ub{{\boldsymbol u}}
\def\vb{{\boldsymbol v}}

\def\xb{{\boldsymbol x}}

\def\zb{{\boldsymbol z}}

\def\Gammab{{\boldsymbol      \Gamma}}           
\def\Deltab{{\boldsymbol      \Delta}}

\def\Lambdab{{\boldsymbol     \Lambda}}

\def\Sigmab{{\boldsymbol      \Sigma}}

\def\Omegab{{\boldsymbol      \Omega}}

\newtheorem{thm}{Theorem}[section]

\newtheorem{lem}[thm]{Lemma}

\newcommand{\vecop}{\mathop{\rm{vec}}}
\newcommand{\trop}{\mathop{\rm{tr}}}
\newcommand{\argmin}{\mathop{\rm argmin}}
\newcommand{\argmax}{\mathop{\rm argmax}}
\newcommand{\norm}[1]{\left\lVert #1 \right\rVert}

\newcommand{\Real}{\mathbb{R}}
\newcommand{\real}{\mathbb{R}}
\newcommand{\St}{\mathcal{S}}

\newcommand{\tr}{\mathop{\rm{tr}}}
\def\NN{N}
\def\vec{{\rm vec}}

\def\opt{{\rm opt}}

\def\1{{\bm 1}}
\def\0{{\bm 0}}

\def\diag{{\rm diag}}

\def\intercal{\top}
\def\i{{[i]}}
\def\ione{{[1]}}
\def\itwo{{[2]}}
\def\iN{{[N]}}
\def\sp{{\rm sp}}

\newcommand{\LeftEqNo}{\let\veqno\@@leqno}

\graphicspath{{figures/}}

\begin{document}

\title{Integrating Multiple Random Sketches for Singular Value Decomposition}

\author{Ting-Li Chen$^*$\!\!\! \and
Dawei D. Chang$^\dag$\!\!\! \and
       Su-Yun Huang$^*$\!\!\! \and
       Hung Chen$^\dag$\!\!\! \and
       Chienyao Lin\thanks{Institute of Statistical Science, Academia Sinica,
         Taipei 115, Taiwan ({\tt tlchen@stat.sinica.edu.tw, syhuang@stat.sinica.edu.tw,
         youyuoims94@gmail.com})} \and
       Weichung Wang\thanks{Institute of Applied Mathematical Sciences,
         National Taiwan University, Taipei 106, Taiwan
         ({\tt davidzan830@gmail.com, hchen@math.ntu.edu.tw},
         send correspondence to {\tt wwang@ntu.edu.tw})}}

\maketitle

\begin{abstract}
The singular value decomposition (SVD) of large-scale
matrices is a key tool in data analytics and scientific computing.
The rapid growth in the size of matrices further increases the need
for developing efficient large-scale SVD algorithms. Randomized SVD
based on one-time sketching has been studied, and its potential has
been demonstrated for computing a low-rank SVD. Instead of exploring different single random sketching techniques, we propose a Monte Carlo type integrated SVD algorithm based on multiple random sketches.
The proposed integration
algorithm takes multiple random sketches and then integrates the results
obtained from the
multiple sketched subspaces. So that the integrated SVD can achieve
higher accuracy and lower stochastic variations.
The main component of the integration
is an optimization problem with a matrix Stiefel manifold
constraint. The optimization problem is solved using
Kolmogorov-Nagumo-type averages. Our theoretical analyses show that
the singular vectors can be induced by population averaging and
ensure the consistencies between the computed and true subspaces and
singular vectors. Statistical analysis further proves a strong Law
of Large Numbers and gives a rate of convergence by the Central
Limit Theorem. Preliminary numerical results suggest that the
proposed integrated SVD algorithm is promising.
\end{abstract}



\section{Introduction}

Singular value decomposition (SVD) of a matrix has been an essential
tool in various theoretical studies and practical applications for
decades. SVD has been studied in the fields of numerical linear
algebra, applied mathematics, statistics, computer sciences, data
analytics, physical sciences, engineering, and others. Its
applications include imaging, medicine, social networks, signal
processing, machine learning, information compression, and finance,
just to name a few. In this article, we focus on low-rank SVD of a
matrix, rather than the full SVD, which is sufficient in many
scenarios. We concern the rank-$k$ SVD of a given ${m\times n}$ real
matrix
\begin{equation}
\Ab =\Ub\Sigmab\Vb^\intercal
 \approx {\Ub}_{k}\, {\Sigmab}_{k}\,{\Vb}_{k}^\intercal,
\label{eq:rank_k_svd}
\end{equation}
where $\Ub\Sigmab\Vb^\intercal$ is the full SVD and ${\Ub}_{k}\,
{\Sigmab}_{k}\,{\Vb}_{k}^\intercal$ is the truncated rank-$k$ SVD.
The columns of ${\Ub}_{k}$ and ${\Vb}_{k}$ are the $k$ leading left
and right singular vectors of $\Ab$, respectively. The diagonal
entries of ${\Sigmab}_{k}$ are the $k$ largest singular values of
$\Ab$.

The role of SVD in all types of applications remains vivid in the
current big data era. However, as the size of data matrices
generated or collected from simulations, experiments, detections,
and observations continues to increase quickly, it is becoming a
challenge to compute SVD for large matrices.

Randomized algorithms have been proposed and studied recently to
find a rank-$k$ SVD of a large matrix. Such algorithms randomly
project or sample from the underlying matrix $\Ab$ to obtain a
reduced matrix in a low-dimensional subspace. SVD of this reduced
matrix in low dimensions is performed and then mapped back to the
original space to obtain an approximate SVD of $\Ab$. Several survey
papers have reviewed the algorithms, computed error bounds, and
numerical experiments in detail
\cite{slide:RandNLA,halko2012randomized,mahoney2011randomized,slide:gunnar,woodruff2014sketching,zhang2015singular}.
Randomized SVD algorithms have been not only investigated
theoretically but also used to solve problems such as large data set
analysis \cite{halko2011algorithm}, overdetermined least squares
\cite{rokhlin2008fast}, partial differential equations
\cite{sabelfeld2012stochastic}, and inverse problems
\cite{xiang2013regularization}. Randomized SVD is also used to solve
linear system problems
\cite{strohmer2006randomized,xia2012superfast} or act as a
preconditioner \cite{coakley2011fast,
demanet2012matrix,grigori2014robust,gu2015subspace,qian2008stochastic}.
Randomized SVD implementations built on top of MATLAB
\cite{szlam2014implementation}, parallel computers
\cite{skylark,halko2012randomized,yang2016implementing}, and
emerging architectures, such as GPUs \cite{yamazaki2015randomized},
have appeared and continue to be improved.

Randomized SVD algorithms use a single random sketch and have
demonstrated their advantages in various situations. In this
article, we enhance randomized-type SVD by proposing a randomized
SVD that integrates results obtained from multiple random sketches
such that the integrated rank-$k$ SVD can achieve higher accuracy
and less stochastic variations. Our discussion and analysis of this
new algorithm include the following.
\begin{itemize}
\item
We introduce a new way to compute rank-$k$ SVD based on multiple
random sketches in Algorithm~\ref{alg:isvd}. Then, we develop
two equivalent optimization problems (Theorem~\ref{lemma:q_ell})
with the Stiefel manifold constraint to integrate multiple
subspace information sources from multiple random sketches. The
proposed algorithm can be viewed as a Monte Carlo method that
randomly samples many subspaces with an integration procedure
based on averaging.

\item
To integrate the results obtained using multiple sketches, we
propose Algorithm~\ref{alg:KN} to solve the constrained
optimization problem iteratively and analyze its convergence
behavior. In each iteration, the algorithm moves the current
iterate using Kolmogorov-Nagumo-type averaging on top of the
Stiefel manifold. The approach is motivated by averaging the
independently identically distributed (i.i.d.) results from
multiple sketches.

\item
The algorithm is analyzed statistically. In a key argument shown
in Theorem~\ref{lemma:key}, we assert that singular vectors can
be induced by population averaging. This key theorem connects
the two subspaces formed by the integration process and by the
true singular vectors. Furthermore, based on this key theorem,
we are able to prove the consistencies in terms of the subspace
and singular vectors, the strong Law of Large Numbers, and the
Central Limit Theorem for convergence rate.

\item
The numerical results, such as Figure~\ref{fig:similarity},
suggest that the integrated SVD can achieve higher accuracy and
less stochastic variations using multiple random sketches.

\end{itemize}

This paper is organized as follows. We introduce the integrated SVD
algorithm in Section~\ref{sec:irsvd}, followed by a detailed
discussion regarding how we integrate the multiple sketched results
in Section~\ref{sec:intQ}. We analyze the algorithm statistically in
Section~\ref{sec:stat}. We present numerical results in
Section~\ref{sec:num}. Finally, we conclude the paper in
Section~\ref{sec:con}.

For the notations, we use lowercase letters or Greek letters for
scalars (e.g., $m$, $n$, and $\tau$), bold face letters for vectors
(e.g., $\xb$ and $\ub$), and bold face uppercase letters or bold
face Greek letters for matrices (e.g., $\Ab$ and $\Omegab$). We use
$\|\cdot\|_\sp$ and $\|\cdot\|_F$ to denote the matrix spectral norm
and the matrix Frobenius norm, respectively.
Table~\ref{tab:notation} summarizes the notations used in this
article. We assume $m\le n$ here; however, all the
algorithms and theoretical results can be applied to other cases.

\begin{table}
\centering
 \begin{tabular}{llllllll}
 \hline
 $m$, $n$ & Row and column dimensions of a matrix $\Ab$ with the assumption $m\le n$\\
 $k$ & Desired rank of approximate SVD\\
 $p$ & Oversampling parameter\\
 $\ell$ & Dimension of randomized sketches, i.e., $\ell=(k+p) \ll n$\\
 $q$ & Exponent of the power method in Step~\ref{step:randproj}
       of Algorithms~\ref{alg:rsvd} and \ref{alg:isvd}\\
 $N$ & Number of random sketches in Algorithm~\ref{alg:isvd} \\
 \hline
 $\Ab=\Ub\Sigmab \Vb^\intercal$ &
              An $m\times n$ matrix and its SVD\\
 $\Ab\approx\Ub_k\Sigmab_k \Vb_k^\intercal$ & Rank-$k$ SVD defined
       in~\eqref{eq:rank_k_svd}\\
 $\Ab \approx \widetilde{\Ub}_k \widetilde{\Sigmab}_k \widetilde{\Vb}_k^\intercal$ &
       Rank-$k$ SVD defined in \eqref{eq:rsvd} and computed by Algorithm~\ref{alg:rsvd} (rSVD)\\
 $\Ab \approx \widehat{\Ub}_k \widehat{\Sigmab}_k \widehat{\Vb}_k^\intercal$ &
       Rank-$k$ SVD defined in \eqref{eq:isvd} and computed by Algorithm~\ref{alg:isvd} (iSVD)\\
 \hline
 $\Omegab$ &  A Gaussian random projection matrix in Algorithm~\ref{alg:rsvd} \\
 $\Omegab_\i$, $i=1:N$ & The $i$th Gaussian random projection matrix in Algorithm~\ref{alg:isvd}\\
 $\Qb_\i$, $i=1:N$ &  The $i$th orthonormal basis of the sketched subspace in Algorithm~\ref{alg:isvd}\\
 $\overline\Qb$ & The integrated orthonormal basis of the sketched subspace\\
 $\overline{\Pb}$ & The average of $\Qb_\i\Qb_\i^\intercal$ defined
     in~\eqref{eq:def_P}\\
     & ($\overline\Pb$ is simply for notation usage.
     It is not for computational purposes.)  \\
 \hline
 $\St_{r,c}$ & Matrix Stiefel manifold $\St_{r,c}:=\left\{\Wb\in\Real^{r\times c}:
  \Wb^\intercal\Wb =\Ib_c \mbox{ and } r\ge c\right\}$ \\
 ${\mathcal T}_{\Qb}\St_{r,c}$ & Tangent space of $\St_{r,c}$ at $\Qb$\\
 $\Qb_{c}$ & The current iterate for computing $\overline\Qb$
    in Algorithms~\ref{alg:KN}\\
 $\Qb_{+}$ & The updated iterate for computing $\overline\Qb$
    in Algorithms~\ref{alg:KN}\\
\hline
 $\Qb$, $\Wb$ & Points located on the matrix Stiefel manifold\\
 $\Xb$ & A point located on a tangent space\\
 $F(\Qb)$ & The objective function defined in \eqref{eq:optQ} for computing $\overline\Qb$ \\
 $\Gb_F(\Qb)$ & Gradient of the objective function $F(\Qb)$ \\
 $\Db_F(\Qb)$ & Projected gradient onto the
   tangent space ${\mathcal T}_{\Qb}\St_{m,\ell}$\\
\hline
 $\varphi_{\Qb}$ & A lifting map to the tangent space in terms of $\Qb$\\
 $\varphi_{\Qb}^{-1}$ & A (specified version of) retraction map to the Stiefel manifold \\
  & ($\varphi_{\Qb}$ and $\varphi_{\Qb}^{-1}$ are associated with the Kolmogorov-Nagumo-type average.) \\
\hline
 \end{tabular}
 \caption{Notations used in this article.}
\label{tab:notation}
\end{table}

\section{Singular Value Decomposition via Multiple Random Sketches}
\label{sec:irsvd}

Randomized algorithms have been proposed to compute an approximate
rank-$k$ SVD for matrices arising in various applications
\cite{halko2012randomized,mahoney2011randomized,woodruff2014sketching,zhang2015singular}.
The main idea of these algorithms is to (i)~randomly project the
matrix to a low-dimensional subspace, (ii)~compute the SVD in this
random subspace, and (iii)~map this subspace SVD back to the
original high-dimensional space. If the random sketch can capture
most of the information regarding the largest $k$ singular values
and singular vectors, these algorithms can obtain satisfactory
approximate rank-$k$ SVDs. We briefly review these randomized SVDs
in Section~\ref{sec:single}.

To improve these randomized SVD algorithms based on a {\it single}
sketch, it is natural to ask how we can find a better random
subspace. Instead of exploring different single random sketching
techniques, we propose a Monte Carlo integration method based on
{\it multiple} random sketches in Section~\ref{sec:mult}. The key
idea is to repeat the process of random sketching multiple times.
The multiple low-dimensional subspaces are then integrated. Based on
this integrated subspace, we compute a rank-$k$ approximate SVD
accordingly. By taking multiple random sketches, the resulting
integrated SVD is expected to have higher accuracy and smaller
stochastic variation. On the other hand, the multiple sketches can
be performed on parallel computers to reduce the execution time.
Furthermore, the aforementioned multiple sketches lead to multiple
low-dimensional random subspaces. We present an optimal
representation of the multiple subspaces in Section~\ref{sec:opt},
which is defined by a constrained optimization problem.

\subsection{Single Random Sketch}
\label{sec:single}

Algorithm~\ref{alg:rsvd} is a common procedure for {\it randomized
SVD} (rSVD)~\cite{halko2011finding,rokhlin2009randomized} used to
compute a rank-$k$ approximate SVD
\begin{equation}
\Ab \approx \widetilde{\Ub}_{k}\, \widetilde{\Sigmab}_{k}\,\widetilde{\Vb}_{k}^\intercal.
\label{eq:rsvd}
\end{equation}
The algorithm includes the following steps.
\begin{description}
\item[Step~\ref{step:genrandmtx}.] The algorithm first generates a
  random matrix $\Omegab \in \mathbb{R}^{n\times \ell}$.

\item [Step~\ref{step:randproj}.]
The random matrix $\Omegab$ is used to map the matrix $\Ab$ to a
low-dimensional subspace by $\Yb =
(\Ab\Ab^\intercal)^{q}\Ab\Omegab$. The $q$th power of the matrix
$\Ab\Ab^\intercal$ is applied here to improve the accuracy for
the slow decay singular values
\cite{halko2011finding,witten2013randomized}.

\item[Step~\ref{step:QR}.]
Compute an orthonormal basis of $\Yb$ by, e.g., QR factorization
or SVD so that the matrix $\Qb$ spans the randomized column
subspace of $\Ab\Omegab$.

\item[Steps~\ref{step:lowdimsvd} and \ref{step:updateU}.]
A smaller scale $\ell \times n$ SVD and a matrix multiplication
are performed to compute the SVD of the projected matrix in the
column space of $\Qb$. These two steps are equivalent to the
operation $\Qb \Qb^\intercal \Ab = \widetilde{\Ub}_{\ell}\,
\widetilde{\Sigmab}_{\ell}\, \widetilde{\Vb}_{\ell}^\intercal$.

\item[Step~\ref{step:extractsvd}.] Because the matrices $\widetilde{\Ub}_{\ell}$,
$\widetilde{\Sigmab}_{\ell}$, and $\widetilde{\Vb}_{\ell}$
contain over-sampled singular vectors and singular values, we
extract the largest rank-$k$ approximate SVD including
$\widetilde{\Ub}_{k}$, $\widetilde{\Sigmab}_{k}$, and
$\widetilde{\Vb}_{k}$ from these matrices.
\end{description}

\begin{algorithm}
  \caption{Randomized SVD with a single sketch (rSVD).}
  \label{alg:rsvd}
  \begin{algorithmic}[1]
    \REQUIRE $\Ab$ (real $m \times n$ matrix), $k$ (desired rank of approximate SVD), $p$ (oversampling parameter),
    $\ell=k+p$ (dimension of the sketched column space),
    $q$ (exponent of the power method)
    \ENSURE Approximate rank-$k$ SVD of $\Ab \approx \widetilde{\Ub}_k\,
        \widetilde{\Sigmab}_k \, \widetilde{\Vb}_k^\intercal$
    \STATE Generate an $n \times \ell$ random matrix $\Omegab$ \label{step:genrandmtx}  
    \STATE Assign $\Yb \leftarrow (\Ab\Ab^\intercal)^{q}\Ab\Omegab$ \label{step:randproj}
    \STATE Compute $\Qb$ whose columns are an orthonormal basis of $\Yb$ \label{step:QR}
    \STATE Compute the SVD of $\Qb^\intercal\Ab = \widetilde{\Wb}_{\ell}\,
           \widetilde{\Sigmab}_{\ell}\, \widetilde{\Vb}_{\ell}^\intercal$ \label{step:lowdimsvd}
    \STATE Assign $\widetilde{\Ub}_{\ell} \leftarrow \Qb\widetilde{\Wb}_{\ell}$ \label{step:updateU}
    \STATE Extract the largest $k$ singular-pairs from $\widetilde{\Ub}_{\ell}$, $\widetilde{\Sigmab}_{\ell}$,
     $\widetilde{\Vb}_{\ell}$ to obtain $\widetilde{\Ub}_k$, $\widetilde{\Sigmab}_k$, $\widetilde{\Vb}_k$ in \eqref{eq:rsvd}
       \label{step:extractsvd}
  \end{algorithmic}
  \end{algorithm}

\subsection{Multiple Random Sketches}
\label{sec:mult}

The rSVD (Algorithm~\ref{alg:rsvd}) maps the matrix $\Ab$ onto a
low-dimensional subspace using a {\it single} random sketch. We
extend the rSVD by proposing an integrated singular value
decomposition (iSVD), which uses {\it multiple} sketches. The
procedure of iSVD is outlined in Algorithm~\ref{alg:isvd}. In
addition to those input parameters listed in rSVD,  the proposed
iSVD (Algorithm~\ref{alg:isvd}) takes an extra parameter: the number
of random sketches $N$. In return, iSVD outputs the integrated
approximate rank-$k$ SVD
\begin{equation}
\Ab \approx \widehat{\Ub}_{k}\, \widehat{\Sigmab}_{k}\,\widehat{\Vb}_{k}^\intercal.
\label{eq:isvd}
\end{equation}

\begin{algorithm}
\caption{Integrated SVD with multiple sketches (iSVD).}
\label{alg:isvd}
\begin{algorithmic}[1]
  \REQUIRE  $\Ab$ (real $m \times n$ matrix), $k$ (desired rank of approximate SVD),
  $p$ (oversampling parameter), $\ell=k+p$
  (dimension of the sketched column space),
    $q$ (exponent of the power method), $N$ (number of random sketches)
  \ENSURE Approximate rank-$k$ SVD of $\Ab \approx
       \widehat{\Ub}_k \widehat{\Sigmab}_k \widehat{\Vb}^\intercal  _k$
  \STATE Generate $n \times \ell$ random matrices $\Omegab_\i$ for $i=1,\ldots, N$ \label{line:isvd_rand}
  \STATE Assign $\Yb_\i \leftarrow (\Ab\Ab^\intercal)^{q}\Ab\Omegab_\i$ for $i=1,...,N$
         (in parallel) \label{line:isvd_proj}
  \STATE Compute $\Qb_\i$ whose columns are an orthonormal basis of $\Yb_\i$ (in parallel) \label{line:isvd_Q}\\
  \STATE Integrate $\overline{\Qb} \leftarrow \{\Qb_\i\}_{i=1}^{N}$
         (by Algorithm~\ref{alg:KN}) \label{step2:intQ}\\
  \STATE Compute the SVD of $\overline{\Qb}^\intercal\Ab =
      \widehat{\Wb}_{\ell}\, \widehat{\Sigmab}_{\ell}\, \widehat{\Vb}_{\ell}^\intercal$ \label{step2:lowdimsvd}
  \STATE Assign $\widehat{\Ub}_{\ell} \leftarrow \overline{\Qb}\widehat{\Wb}_{\ell}$ \label{step2:updateU}
  \STATE Extract the largest $k$ singular pairs from $\widehat{\Ub}_{\ell}$, $\widehat{\Sigmab}_{\ell}$, $\widehat{\Vb}_{\ell}$
      to obtain $\widehat{\Ub}_k$, $\widehat{\Sigmab}_k$, $\widehat{\Vb}_k$ in \eqref{eq:isvd} \label{line:isvd_result}
\end{algorithmic}
\end{algorithm}

In Steps~\ref{line:isvd_rand}, \ref{line:isvd_proj}, and
\ref{line:isvd_Q} of Algorithm~\ref{alg:isvd}, iSVD performs
multiple sketches by repeating the sketching process described in
the first three steps of the rSVD algorithm $N$ times. In
Step~\ref{step2:intQ}, the multiple orthonormal basis matrices
$\Qb_\i$ are integrated. Using the integrated orthonormal basis
matrix $\overline{\Qb}$, we can obtain an approximate SVD by
Steps~\ref{step2:lowdimsvd}, \ref{step2:updateU},
and~\ref{line:isvd_result}. Note that, in Step~\ref{step:genrandmtx}
of the two algorithms, we consider Gaussian random projection
matrices $\Omegab$ in rSVD and $\Omegab_\i$ in iSVD. Either
$\Omegab$ or $\Omegab_\i$ is an $n\times \ell$ random matrix whereby
each of the entries is i.i.d. standard Gaussian. The matrix $\Ab
\Omegab$ (or $\Ab \Omegab_\i$) is a random mapping from
$\mathbb{R}^{m\times n}$ to a low-dimensional subspace $\Yb
(\mbox{or}\ \Yb_\i) \in \mathbb{R}^{m\times \ell}$ with \mbox{$\ell
\ll n$}. Furthermore, each of the columns in $\Ab \Omegab$ (or $\Ab
\Omegab_\i$) is a linear combination of the columns of $\Ab$ with
random Gaussian mixing coefficients.

\begin{figure}
\begin{center}
\subfloat[$\Qb_\i\in \real^{3\times1}$ ($k=\ell=1$) and $N=2$.]{\includegraphics[width=.42\textwidth]{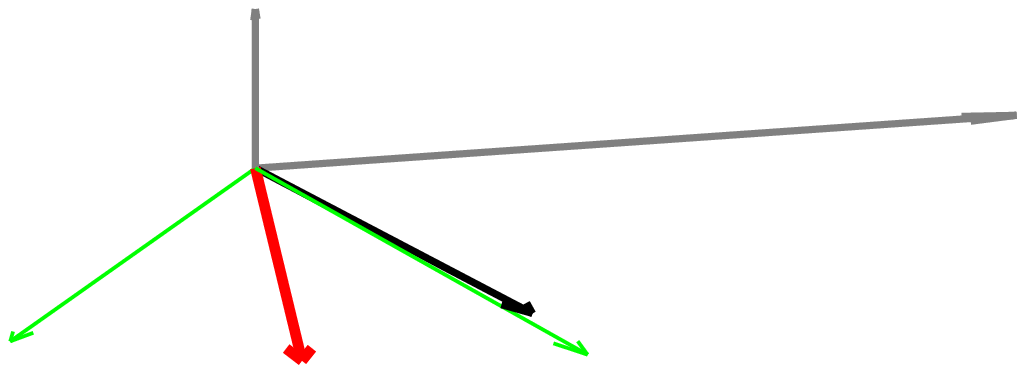}\label{subfig:isvd_1_2}}\ \ \ \ \ \
\subfloat[$\Qb_\i\in \real^{3\times1}$ ($k=\ell=1$) and $N=5$.]{\includegraphics[width=.42\textwidth]{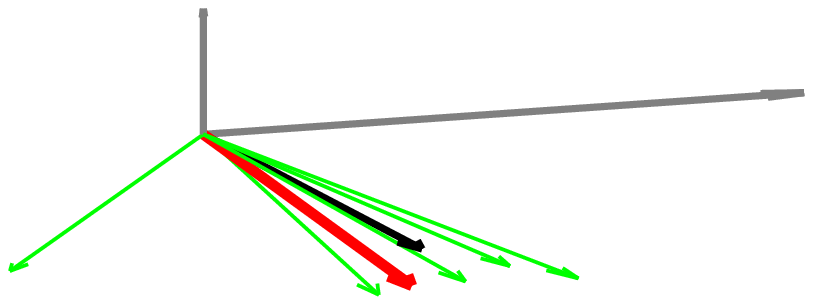}\label{subfig:isvd_1_5}} \\
\subfloat[$\Qb_\i\in \real^{3\times2}$ ($k=\ell=2$) and $N=2$.]{\includegraphics[width=.4\textwidth]{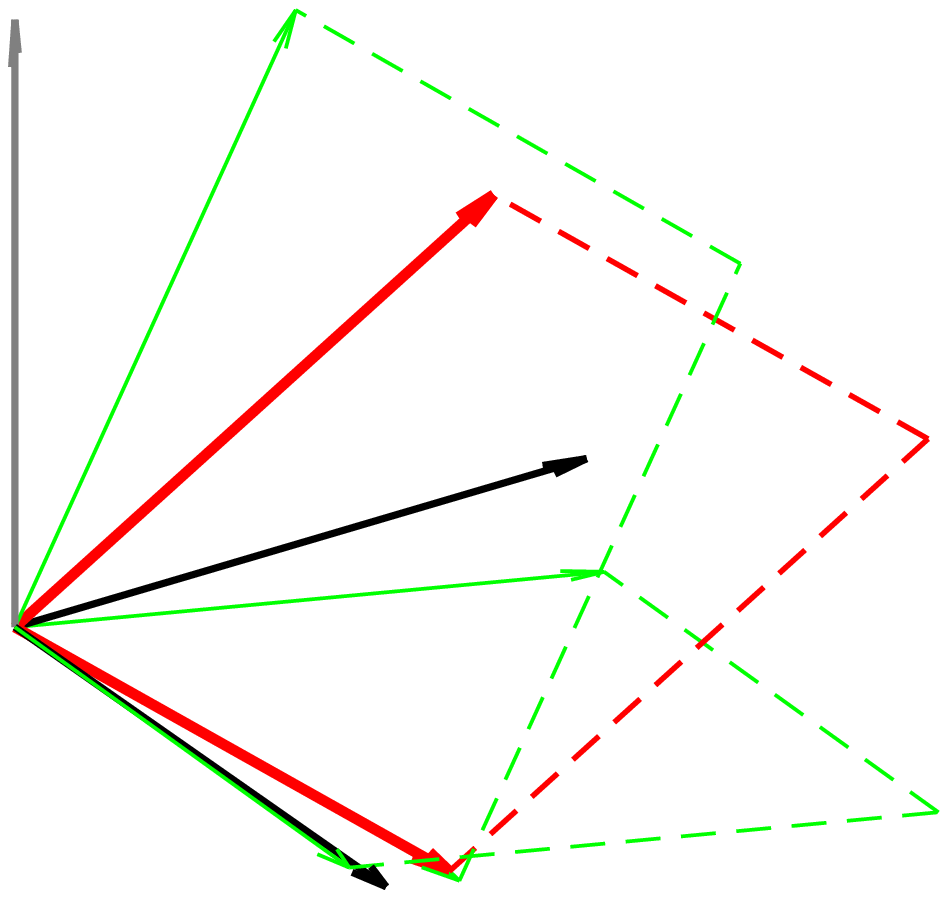}\label{subfig:isvd_2_2}}\ \ \ \ \ \
\subfloat[$\Qb_\i\in \real^{3\times2}$ ($k=\ell=2$) and $N=5$.]{\includegraphics[width=.4\textwidth]{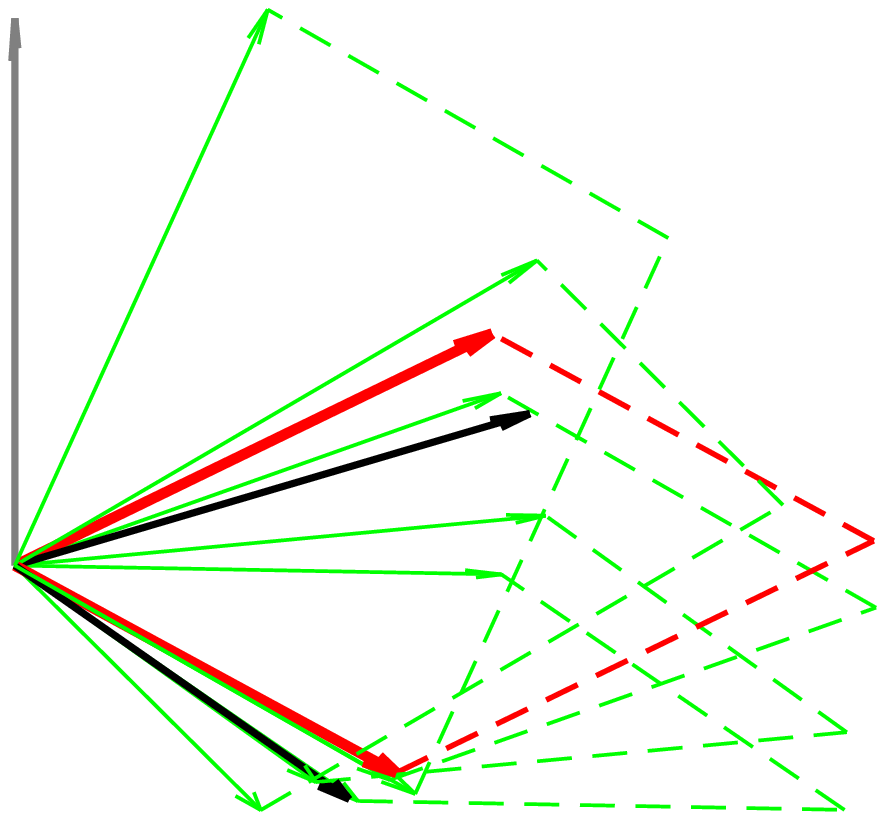}\label{subfig:isvd_2_5}}
\end{center}
\caption{The orthonormal bases of randomly projected $1$-dimensional
subspaces $\Qb_\i$, $i=1,\dots,N$, are plotted (in green) for (a) $N=2$ and (b) $N=5$.
The basis of the integrated subspace $\overline{\Qb}$
(in red) is closer to the first singular vector $[1,0,0]$ for $N=5$.
Similar plots on $2$-dimensional subspaces are plotted for (c) $N=2$
and (d) $N=5$. The $2$-dimensional subspace spanned by
$\overline{\Qb}$ is closer to the subspace spanned by the first two
singular vectors for $N=5$.} \label{fig:isvd}
\end{figure}

We use a simple example to illustrate the ideas of iSVD. Let
$\Ab=\diag([25, 5, 1])$ be a diagonal $3\times 3$ matrix. We have
the (true) SVD of $\Ab=\Ub \Sigmab \Ub^\intercal,$ where $\Ub$ is
the $3\times 3$ identify matrix and $\Sigmab=\diag([25, 5, 1])$. As
shown in Parts~(a) and~(b) of Figure~\ref{fig:isvd}, we randomly
project $\Ab$ onto $\real^{3\times1}$ by letting $\Omegab_\i \in
\real^{3\times1}$ with $q=0$ and $N=2, 5$. The bases of the
projected subspaces $\Qb_\i\in \real^{3\times1}$ are shown by the
green vectors, and the integrated basis $\overline{\Qb}$ is shown as
the red vector. It is clear that the $\overline{\Qb}$ corresponding
to $N=5$ is more close to the first singular vector $[1,0,0]$.
Consequently, we can obtain more accurate SVD $\widehat{\Ub}_{k}\,
\widehat{\Sigmab}_{k}\,\widehat{\Vb}_{k}^\intercal$ by iSVD over the
subspace $\overline{\Qb}$. In Parts~(c) and~(d), we show  similar
results obtained by letting $\Omegab_\i \in \real^{3\times2}$. The
$\overline{\Qb}$ is more close to the $2$-dimensional subspace
spanned by the first two singular vectors $[1,0,0]$ and $[0,1,0]$
using larger $N$.

We have proposed the iSVD algorithm based on multiple random
sketches in Section~\ref{sec:mult}. Obviously, the key component of
 iSVD is the integration process in Step~\ref{step2:intQ} of
Algorithm~\ref{alg:isvd}. This is the focus of the next section.

\subsection{An Optimal Representation of the Multiple Projected Subspaces}
\label{sec:opt}

The integration process in Step~\ref{step2:intQ} of
Algorithm~\ref{alg:isvd} finds a matrix $\overline{\Qb}$ that
``best'' represents the matrices $\Qb_\i$ for $i=1,\ldots,N$. In
other words, because each $\Qb_\i$ contains the orthonormal basis of
the randomly projected subspace $\Yb_\i$, the process intends to
integrate these $N$ randomly projected subspaces into a single
subspace spanned by the columns of $\overline{\Qb}$. Consequently,
this integrated subspace contains as much information of the leading
left singular vectors as possible. Then, in
Steps~\ref{step2:lowdimsvd} and \ref{step2:updateU} of
Algorithm~\ref{alg:isvd}, we compute the SVD in $\overline{\Qb}\,
\overline{\Qb}^\intercal\Ab$, which is the low-dimensional
projection of $\Ab$ onto the subspace spanned by the columns of $\overline{\Qb}$. In
particular, we define such best representation $\overline{\Qb}$ by
solving the following optimization problem:
\begin{equation}
\overline{\Qb} :=\argmin_{\Qb\in \St_{m,\ell}}\sum_{i=1}^N
\left\|\Qb_\i\Qb_\i^\intercal -\Qb\Qb^\intercal\right\|_F^2.
\label{eq:QiQ_1}
\end{equation}
The matrix $\overline{\Qb}$ is constrained on the matrix Stiefel
manifold because the columns of $\overline{\Qb}$ form the
orthonormal basis of the integrated subspace. Next, we justify this
definition of the $\overline{\Qb}$ from the viewpoints of geometry
and stochastic expectation.

This definition of $\overline{\Qb}$ has its geometrical motivations.
At first glance, we can average the $\Qb_\i$ by computing $\Qb_{ave}
= \frac{1}{N}\sum_{i = 1}^{\NN} \Qb_\i$ and then orthonormalize the
columns of $\Qb_{ave}$ to obtain an average representation of
$\Qb_\i$. However, this simple averaging scheme can be misguided.
For example, let $\Qb_\ione = [\qb_{1},\qb_{2}]$ and $\Qb_\itwo =
[\qb_{2},\qb_{1}]$ be two equivalent matrices with respect to an
orthogonal transformation. The simple averaging schemes suggests
that $\Qb_{ave} = \frac{1}{2}[\qb_{1}+\qb_{2}, \qb_{1}+\qb_{2}]$,
which is rank deficient. Fortunately, we observe that the
equivalence of $\Qb_\ione$ and $\Qb_\itwo$ can be revealed by the
equation $\Qb_\ione \Qb_\ione^\intercal = \Qb_\itwo
\Qb_\itwo^\intercal$. On the other hand, in general, any orthogonal
transformation on the right-hand side of $\Qb_\i$ can be represented
by $\Qb_\i\Rb$, where $\Rb$ is an orthogonal matrix. The fact that
$(\Qb_\i\Rb) (\Qb_\i\Rb)^\intercal = \Qb_\i \Qb_\i^\intercal$
suggests that the matrices $\Qb_\i$ and $\Qb_\i \Rb$ are equivalent
in the sense that they span the same column subspace. These
geometric observations partially motivate us to define the
integrated orthonormal matrix $\overline{\Qb}$ shown in
\eqref{eq:QiQ_1}.

Furthermore, we emphasize another important reason why we focus on
the projection matrices $\Qb_\i\Qb_\i^\intercal$ by presenting the
following key Theorem~\ref{lemma:key}. The theorem suggests that the
population average of the projection matrices
$E(\Qb_\i\Qb_\i^\intercal)$ can reveal the {\it true} left singular
vectors of $\Ab$. Furthermore, in Section~\ref{sec:stat}, we will
apply Theorem~\ref{lemma:key} to show the Strong Law of Large
Numbers, the consistency of the singular vectors, and the Central
Limit Theorem for iSVD. The proof of Theorem~\ref{lemma:key} can be found in Appendix~\ref{apdnx:lemma:key}.
\begin{thm}[Singular vectors induced by population averaging]\label{lemma:key}
Let $\Qb_{[i]}$ be the orthonormal basis of the $i$th random
subspace computed by Algorithm~\ref{alg:isvd} with $\Omegab_\i$
having i.i.d. Gaussian entries. At the population level, the
expected arithmetic mean of these projection matrices has the
property of having the same left singular vectors as the matrix
$\Ab$. Specifically,
\begin{equation}
E (\Qb_{[i]}\Qb_{[i]}^\intercal ) = \Ub\Lambdab\Ub^\intercal.
\label{eq:EQQT}
\end{equation}
Here, $\Ub$ consists of left singular vectors of $\Ab$, as shown
in~\eqref{eq:rank_k_svd}. For $\Lambdab$ defined in
\eqref{eq:Lambda}, we can show that (a)~$\Lambdab$ is a diagonal
matrix, (b)~each of the diagonal entries belongs to $(0,1)$, and
(c)~these diagonal entries are strictly decreasing if the underlying
matrix $\Ab$ has strictly decreasing singular values.
\end{thm}

It is worth noting that we can regard $E(\Qb_{[i]}\Qb_{[i]}^\intercal)$ as the limiting case of taking an average over infinitely many projection matrices: $E(\Qb_{[i]}\Qb_{[i]}^\intercal)=\lim_{N\to \infty}
\overline{\Pb}$, where
\begin{equation}
\overline{\Pb} := \frac 1N  \sum_{i=1}^N \Qb_{[i]}\Qb_{[i]}^\intercal
\label{def:Pbar}
\end{equation}
is the empirical arithmetic average of the $N$ projection matrices.
The theorem suggests the following essential
property. Even a projection matrix $\Qb_\i\Qb_\i^\intercal$ is
associated with a low-dimensional (rank-$\ell$) subspace only, the
arithmetic average $\overline{\Pb}$ contains not only information of
the leading rank-$\ell$ subspace but also information for other
subspaces spanned by {\it all} {\it true} singular vectors if the
number of random sketches $N$ is sufficiently large. Furthermore,
because the entries of $\Lambdab$ are strictly decreasing, the
columns of $\Ub$ match the left singular vectors in a correct order.
The following example illustrates the property given in
Theorem~\ref{lemma:key}. Let $\Ab=\diag([25, 5, 1 , 0.2])$ be a
diagonal $4\times 4$ matrix. Then, we have the true SVD of $\Ab=\Ub
\Sigmab \Ub^\intercal,$ where $\Ub$ is the $4\times 4$ identify
matrix and $\Sigmab=\diag([25, 5, 1  , 0.2])$. By computing the SVD
of $\overline{\Pb} \in \real^{4\times 4}$ with a set of $\Qb_{[i]}\in
\real^{4\times 2}$ (i.e., the dimension of the random sketches
$\ell=2$), we obtain the following estimations of
$$
\small 
\Ub \approx
\left[
\begin{array}{rrrr}
-1.00  &  .021 &  -.004 &   .010 \\
-.021  & -.998 &  -.044 &  -.046 \\
 .004  &  .046 &  -.998 &  -.036 \\
 .009  & -.044 &  -.038 &   .998
\end{array} \right] \mbox{ and }
\left[
\begin{array}{rrrr}
-1.00 & -.020 &  .002 &  .004 \\
 .020 & -1.00 & -.003 & -.009 \\
-.002 &  .003 & -1.00 & -.001 \\
 .004 & -.009 & -.001 &  1.00
\end{array} \right],
$$
for $N=10$ and $N=100$, respectively. The corresponding estimations
of $\Lambdab \approx \diag([.997,\ .848,\ .149,\ .006])$ and
$\diag([.992,\ .819,\ .177,\ .012])$. The
approximate $\Ub$ for $N=100$ is much closer to the whole true $\Ub$
even though the dimension of the random sketches is 2, rather than
the dimension of the matrix, which is 4.

Although Equation~\eqref{eq:EQQT} reveals  important insight into
the average of the $\Qb_\i\Qb_\i^\intercal$, the equation has its
limits from the perspective of numerical computation. First, we
cannot compute the true singular values $\Sigmab$ and the right
singular vectors $\Vb$ using Equation~\eqref{eq:EQQT}. Second, the
diagonal entries of $\Lambdab$ are clustered in the interval
$(0,1)$, and such clustering may affect the accuracy of the computed
$\Ub$. These difficulties can be overcome by considering another
optimization problem shown in Theorem~\ref{lemma:q_ell}.

In Theorem~\ref{lemma:q_ell}, we present two alternative
optimization problems that are equivalent to the optimization
problem \eqref{eq:QiQ_1}. The first equivalent optimization problem
is shown in \eqref{eq:def_P}. In this formulation, we apply the
concept of \eqref{eq:EQQT} to compute $\overline{\Qb}$. The second
equivalent formulation is shown in \eqref{eq:optQ}, which is defined
by a differentiable objective function. We will develop an algorithm
to solve this problem in Section~\ref{sec:intQ}. This decision is
based on the following two reasons. The dimension of
$\Qb_\i\Qb_\i^\intercal$ in \eqref{eq:def_P} ($m \times m$)
can be much larger then the dimension of the matrix
$\Qb^\intercal\overline{\Pb}\Qb$ in \eqref{eq:optQ} ($\ell
\times \ell$). Furthermore, the objective function $F(\Qb)$ is
differentiable, which allows us to develop algorithms to solve the
optimization problem based on the gradient of $F(\Qb)$. The proof of
Theorem~\ref{lemma:q_ell} can be found in
Appendix~\ref{a.proof.equiv}.
\begin{thm}[Equivalent optimization problems] \label{lemma:q_ell}
The minimization problem~(\ref{eq:QiQ_1}) is equivalent to the
following two optimization problems:
\begin{equation}\label{eq:def_P}
(i)\
\overline{\Qb} :=
\argmin_{\Qb\in \St_{m,\ell}}  \left\|\overline{\Pb}
 -\Qb\Qb^\intercal\right\|_F^2
\end{equation}
and
\begin{equation}\label{eq:optQ}
(ii)\
\overline{\Qb}:=\argmax_{\Qb\in \St_{m,\ell}} F(\Qb),
\end{equation}
 where $\overline{\Pb}$ is defined in \eqref{def:Pbar} and $F(\Qb) = \frac12\tr\left(\Qb^\intercal\overline{\Pb}\Qb\right)$ is a
differential function. These equivalences are up to an orthogonal
matrix multiplied on the right-hand side of $\overline{\Qb}$.
\end{thm}

In short, in Step~\ref{step2:intQ} of Algorithm~\ref{alg:isvd}, we
need to integrate the orthogonal matrices $\Qb_\i$ into one
orthogonal matrix $\overline{\Qb} \in \mathbb{R}^{m\times \ell}$
that ``best" represents these matrices. We propose using the
particular $\overline{\Qb}$ defined by the
Stiefel-manifold-constrained optimization problem \eqref{eq:optQ}.
In the next section, we will discuss how we compute $\overline{\Qb}$
by an iterative method based on the Kolmogorov-Nagumo-type average.
In Section~\ref{sec:stat}, we will prove that this optimal
representation $\overline{\Qb}$ converges to the best rank-$\ell$
approximation with probability one when the number of sketches $N$
tends to infinity.

\section{Integration of Sketched Subspaces}
\label{sec:intQ}

The goal of this section is to develop an iterative method based on
a Kolmogorov-Nagumo-type average to compute $\overline{\Qb}$ by
solving the constrained optimization \eqref{eq:optQ}. We start the
development of the integration algorithm by introducing some
background in Section~\ref{sec:kn_update}. This background includes
the Kolmogorov-Nagumo-type average of sample points on a matrix
Stiefel manifold and derives the gradients of the objective
function. Because a Kolmogorov-Nagumo-type average is defined by a
lifting map and a corresponding retraction map, we derive a
particular lifting and retraction map pair that can be applied to
solve the constrained optimization \eqref{eq:optQ} in
Sections~\ref{sec:lifting} and \ref{sec:retraction}. Based on the
lifting and retraction maps, we propose the integration algorithm in
Section~\ref{sec:int_alg}. The convergence analysis and some remarks
on the algorithm are given in Section~\ref{sec:isvd_conv}.

\subsection{Background}
\label{sec:kn_update}

We introduce the Kolmogorov-Nagumo-type average and derive the
gradient of the objective function that will be used to
integrate the sketched subspaces by solving the optimization problem
\eqref{eq:optQ}.

First, we introduce the Kolmogorov-Nagumo-type average. Taking an
average is the most commonly used summary statistic for
independently and identically distributed (i.i.d.) data. The
orthogonal matrices $\{\Qb_{[i]}\}_{i=1}^N$ from repeated runs of
random sketches are independently obtained from a common stochastic
randomization mechanism and thus are i.i.d. The integration of
multiple random sketches can be seen as an ``average'' of these
orthogonal matrices. However, it is no longer in the traditional
sense of taking an average in a Euclidean space; rather, it is a
Kolmogorov-Nagumo-type average defined in~\eqref{KN_ave}. A
Kolmogorov-Nagumo-type average of $\{\mub_i\}_{i=1}^N$ is defined as
\begin{equation}\label{KN_ave}
{\overline \mub}_{KN} =\varphi^{-1}\left( \frac 1 N\sum_{i=1}^N
 \varphi(\mub_i)\right),
\end{equation}
where $\varphi$ is a continuous and locally one-to-one lifting map
and $\varphi^{-1}$ is the paired retraction map. Note that the
traditional arithmetic average can be defined by letting $\varphi$
and $\varphi^{-1}$ be the identity maps, and the geometric average
of positive numbers can be defined by letting $\varphi$ be the
logarithm function and $\varphi^{-1}$ be the exponential function.
In our integration algorithm, we consider the case in which $\mub_i
= \Qb_\i$ and use the notation $\varphi_\Qb$ and
$\varphi_{\Qb}^{-1}$ to emphasize that the lifting and retraction
maps depend on a given $\Qb\in\St_{m,\ell}$. In the next two
sections, we derive a lifting map
$\varphi_\Qb:\St_{m,\ell}\rightarrow{\cal T}_\Qb\St_{m,\ell}$ and
its corresponding retraction map $\varphi_{\Qb}^{-1}: {\cal
T}_\Qb\St_{m,\ell} \rightarrow \St_{m,\ell}$, which satisfy certain
properties for solving the optimization problem \eqref{eq:optQ}.
Various Kolmogorov-Nagumo-type averages of sample points on a
Stiefel manifold can be found in \cite{fiori2014KN,kaneko2013KN}.

Second, we address the (projected) gradient of $F(\Qb)$. Many
optimization schemes, including the scheme to be proposed in
Section~\ref{sec:intQ}, require the derivatives of the objective
function. Theorem~\ref{lemma:q_ell} has asserted that
\eqref{eq:QiQ_1} is equivalent to the  problem \eqref{eq:optQ} with
a differentiable objective function $F(\Qb)$. We further present
Theorem~\ref{thm:proj_grad} to show how we can compute the gradient
ascent direction of $F(\Qb)$ at a certain $\Qb \in
\real^{m\times\ell}$ by~\eqref{def:GF} and show how we can project
the gradient ascent direction to the tangent space of $\St_{m,
\ell}$ at $\Qb$ (denoted as ${\mathcal T}_{\Qb}\St_{m, \ell}$), as
shown in \eqref{def:DF01}. See Appendix~\ref{a.proof.proj_grad} for
the proof of Theorem~\ref{thm:proj_grad}.
\begin{thm} \label{thm:proj_grad}
Let $\Gb_F(\Qb)$ denote the gradient (the usual derivative in the
Euclidean space) of $F(\Qb)$ with respect to
$\Qb\in\real^{m\times\ell}$, and let $\Db_F(\Qb)$ denote the
projected gradient of $\Gb_F(\Qb)$ onto the tangent space ${\mathcal
T}_{\Qb}\St_{m, \ell}$. We have
\begin{equation}
\Gb_F(\Qb) := \left[ \frac{\partial F(\Qb)}{\partial \Qb} \right]
 = \overline{\Pb}\Qb \in \mathbb{R}^{m\times \ell},
\label{def:GF}
\end{equation}
where $\overline{\Pb}$ is defined in \eqref{def:Pbar}, and
\begin{equation}
\Db_F(\Qb) := \Pi_{{\cal T}_\Qb} \Gb_F(\Qb)
  = (\Ib_m  - \Qb\Qb^\intercal )\Gb_F(\Qb),  \label{def:DF01}
\end{equation}
where $\Pi_{{\cal T}_\Qb}$ is the projection from $\real^{m\times
\ell}$ to ${\cal T}_\Qb\St_{m,\ell}$.
\end{thm}

\subsection{The Lifting Map}
\label{sec:lifting}

For a given $\Qb\in\St_{m,\ell}$, we define the lifting map
\begin{equation} \label{def:varphi_Q}
\varphi_{\Qb} (\Wb)  = (\Ib_m-\Qb\,\Qb^\intercal)
\Wb\Wb^\intercal\Qb:\St_{m,\ell}\rightarrow{\cal T}_\Qb\St_{m,\ell}
\end{equation}
for any $\Wb\in\St_{m,\ell}$. Our definition of $\varphi$ leads to
an important property:
\begin{equation}\label{eq:KN_on_tangent}
\frac 1N\sum_{i=1}^N \varphi_{\Qb_c}(\Qb_\i ) = \Db_F(\Qb_c),
\end{equation}
where $\Qb_c$ denotes the current iterate. Specifically, the average
of the mapped points $\varphi_{\Qb_c}(\Qb_\i)$ on the tangent space
of the current iterate is simply the projected gradient at this
current iterate. This property links the Kolmogorov-Nagumo-type
average to the gradient ascent method for the optimal representation
in~(\ref{eq:QiQ_1}) and its equivalent formulation
in~(\ref{eq:optQ}). If the projected gradient $\Db_F(\Qb_c)$ is zero
(or numerically close to zero), then $\Qb_c$ has reached a
stationary point for the optimization problem~(\ref{eq:optQ}). If it
is not zero, we search for the next iterate along the path ${\cal
Q}(\tau):=\varphi^{-1}_{\Qb_c}(\tau \Db_F(\Qb_c))$ on the manifold,
where $\tau>0$ is a step size. In the Kolmogorov-Nagumo-type average, we take $\tau=1$ for simplicity.
Because $\varphi_{\Qb_c}$ is only locally one to one, we need
to specify a version of the retraction map $\varphi_{\Qb_c}^{-1}$ to
pull points on ${\cal T}_{\Qb_c}\St_{m,\ell}$ back to
$\St_{m,\ell}$. Below, we discuss the derivation of a proper version
of $\varphi_{\Qb_c}^{-1}$.

\subsection{The Retraction Map}
\label{sec:retraction}

Next, we derive the corresponding retraction map
$\varphi_{\Qb}^{-1}: {\cal T}_\Qb\St_{m,\ell} \rightarrow
\St_{m,\ell}$. For $\Wb\in \St_{m,\ell}$, we can express it as
$\Wb=\Qb\Cb+\Qb_\bot\Bb$, where $\Qb_\bot^\intercal\Qb_\bot
=\Ib_{m-\ell}$ and $\Qb^\intercal\Qb_\bot=\0$. Without loss of
generality, we may assume that $\Cb$ is symmetric. If not, we can
find an orthogonal matrix $\Rb$ such that $\Cb\Rb$ is symmetric.
Because $(\Wb\Rb)(\Wb\Rb)^\intercal=\Wb\Wb^\intercal$ for any
orthogonal matrix $\Rb$, we treat $\Wb$ and $\Wb\Rb$ as equivalent.
Next, we present two lemmas that will be used in deriving
$\varphi_{\Qb}^{-1}(\Xb)$, where $\Xb\in{\cal T}_{\Qb}\St_{m,\ell}$.
The proofs are given in Appendices~\ref{a.proof.lemma:XX<1/4} and
\ref{a.proof.lemma:C}.
\begin{lem}\label{lemma:XX<1/4}
For a given $\Qb\in\St_{m,\ell}$, we have the following properties.
(a)~The matrix $\frac{\Ib_\ell}4-\varphi_\Qb(\Wb)^\intercal
\varphi_\Qb(\Wb)$ is non-negative definite for any arbitrary
$\Wb\in\St_{m,\ell}$. (b)~Let $\Xb=\frac 1N\sum_{i=1}^N
\varphi_\Qb(\Qb_\i)$. Then, $\frac{\Ib_\ell}4-\Xb^\intercal\Xb$ is
non-negative definite.\\
\end{lem}

\begin{lem}\label{lemma:C}
For a given $\Xb\in{\mathcal T}_{\Qb}\St_{m,\ell}$ that satisfies
the conditions $\Qb^\intercal \Xb = \0$ and $\frac{\Ib_\ell}4
-\Xb^\intercal\Xb$ being non-negative definite, there exists a
$\Wb\in \St_{m,\ell}$ such that $\varphi_\Qb(\Wb)=\Xb$. Furthermore,
if $\Wb$ is restricted to the column span of $\Qb$ and $\Xb$, i.e.,
\[\Wb\in \{\Wb\in\St_{m,\ell}:\Wb=\Qb\Cb+\Xb\Bb,~
\mbox{$\Cb$ symmetric} \},\] then, up to an orthogonal
transformation on the right side, $\Wb$ has to take the following
form $\Wb=\Qb \Cb+ \Xb \Cb^{-1}$, where $\Cb=
\left\{\frac{\Ib_\ell}2 + \left(\frac{\Ib_\ell}4
-\Xb^\intercal\Xb\right)^{1/2}\right\}^{1/2}$ and the matrix square
root is defined in Appendix~\ref{a.matrix_sqrt}.
\end{lem}

Because the matrix square root is not unique, Lemma~\ref{lemma:C}
presents many possible choices of $\Wb$ as a pre-image for $\Xb$
such that $\varphi_\Qb(\Wb)=\Xb$. Here, we will confine the matrix
square root involved in $\Cb$ to be symmetric and non-negative
definite so that the inverse map $\varphi_{\Qb}^{-1}(\Xb) =\Qb \Cb+
\Xb \Cb^{-1}$ is uniquely specified. Furthermore, if
$\Xb\in{\mathcal T}_{\Qb}\St_{m,\ell}$ satisfies the condition that
$\frac{\Ib_\ell}4 -\Xb^\intercal\Xb$ is non-negative definite, then
$\frac{\Ib_\ell}4 -(\tau \Xb^\intercal)(\tau\Xb)$ is also
non-negative definite for any $\tau\in[0,1]$. Thus, we can extend
Lemma~\ref{lemma:C} to obtain a unique path on the manifold, wherein
all matrix square roots involved are taken to be symmetric and
non-negative definite.

\begin{thm}[Retraction Map]\label{thm:C} For a given $\Xb\in{\mathcal
T}_{\Qb}\St_{m,\ell}$ that satisfies the conditions $\Qb^\intercal
\Xb = \0$ and $\frac{\Ib_\ell}4 -\Xb^\intercal\Xb$ being
non-negative definite, there exists a path ${\cal Q}(\tau) \in
\St_{m,\ell}$ for $\tau\in[0,1]$ such that $\varphi_\Qb({\cal
Q}(\tau))=\tau\Xb$ and the retraction map is given by
\begin{equation}\label{eq:phi_inv}
\varphi_\Qb^{-1}(\tau\Xb)=\Qb \Cb+ \tau\Xb \Cb^{-1},
\end{equation}
where
\begin{equation}\label{eq:matrixC}
\Cb=\left\{\frac{\Ib_\ell}2 +
\Big(\frac{\Ib_\ell}4-\tau^2\Xb^\intercal\Xb\Big)^{1/2}\right\}^{1/2}
\end{equation}
with all matrix square roots taken to be symmetric and non-negative
definite.
\end{thm}

\subsection{The Integration Algorithm}
\label{sec:int_alg}

Now, we are ready to propose Algorithm~\ref{alg:KN}, which solves
the optimization problem \eqref{eq:optQ} to find $\overline{\Qb}$ by
iteratively updating the Kolmogorov-Nagumo-type averages. The inputs
of Algorithm~\ref{alg:KN} are the matrices $\{\Qb_\i\}_{i=1}^N \in
\St_{m,\ell}$ and an initial iterate $\Qb_{\rm ini}$. The output of
the algorithm is the (approximate) integrated $\overline\Qb$ defined
in \eqref{eq:optQ}.

\begin{algorithm}
  \caption{Integration of $\{\Qb_\i\}_{i=1}^N$ based on the Kolmogorov-Nagumo-type average.}
  \label{alg:KN}
  \begin{algorithmic}[1]
    \REQUIRE $\Qb_\ione$, $\Qb_\itwo$, $\ldots$, $\Qb_\iN$, $\Qb_{\rm ini}$
    \ENSURE Integrated $\overline\Qb$ defined in \eqref{eq:optQ}
    \STATE Initialize the current iterate $\Qb_{c} \leftarrow \Qb_{\rm ini}$
    \WHILE{(not convergent)} \label{ln:KN_conv}
      \STATE Compute $\overline{\varphi_{\Qb_c}(\Qb_\i)}$ 
        defined in \eqref{eq:ave_Qi} by lifting and averaging\label{step:kn_ave}
      \STATE Perform the retraction mapping $\varphi_{\Qb_c}^{-1}
       \left( \overline{\varphi_{\Qb_c}(\Qb_\i)}\right)$
       to obtain $\Qb_{+}$
       defined in \eqref{def:KN_Q_next} \label{step:kn_pullback}
      \STATE Assign $\Qb_{c} \leftarrow \Qb_{+}$
    \ENDWHILE
  \STATE Output $\overline\Qb = \Qb_c$
  \end{algorithmic}
\end{algorithm}

\begin{figure}
  \begin{center}
  \subfloat[ ]{\includegraphics[width=.55\textwidth]{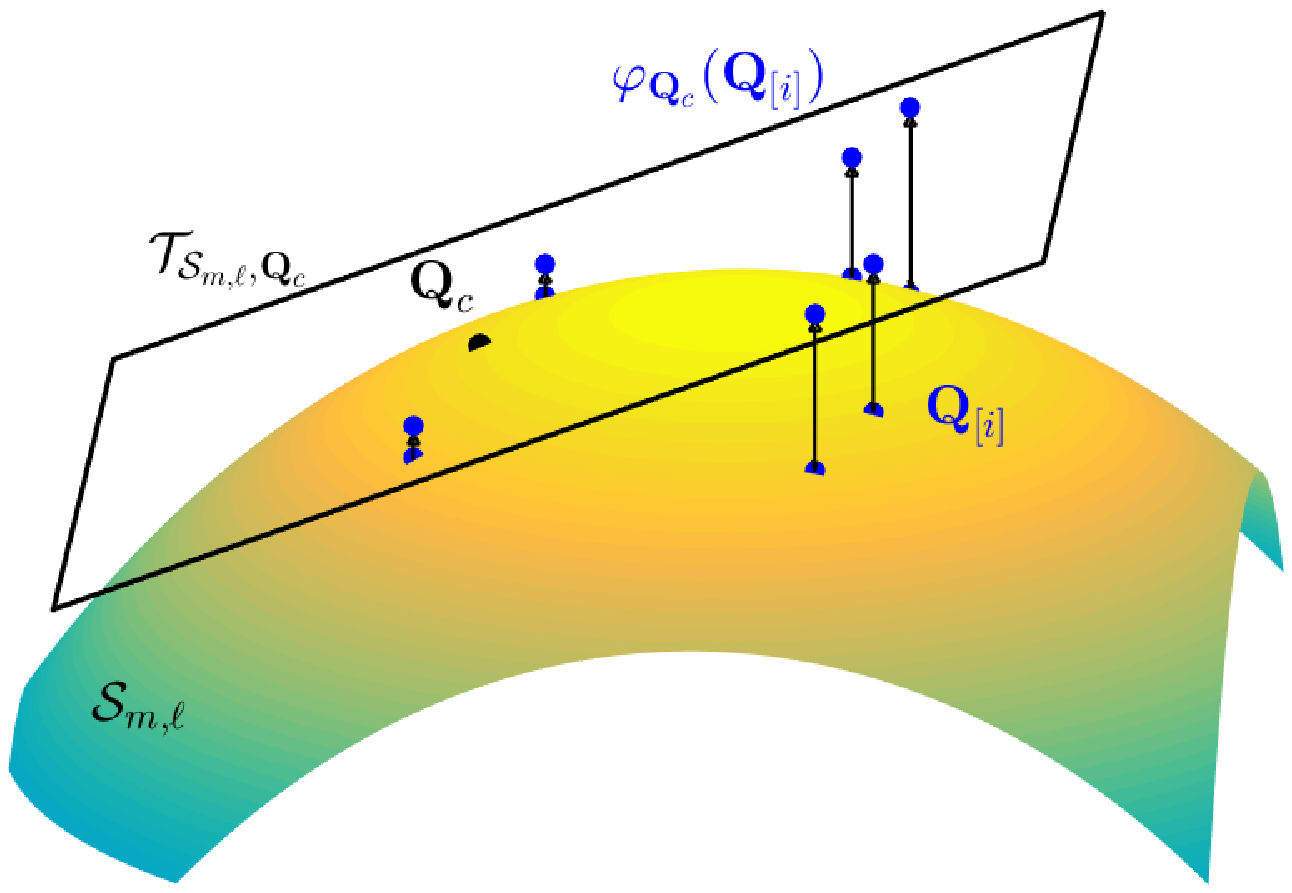}\label{subfig:kn_a}}
  \subfloat[ ]{\includegraphics[width=.55\textwidth]{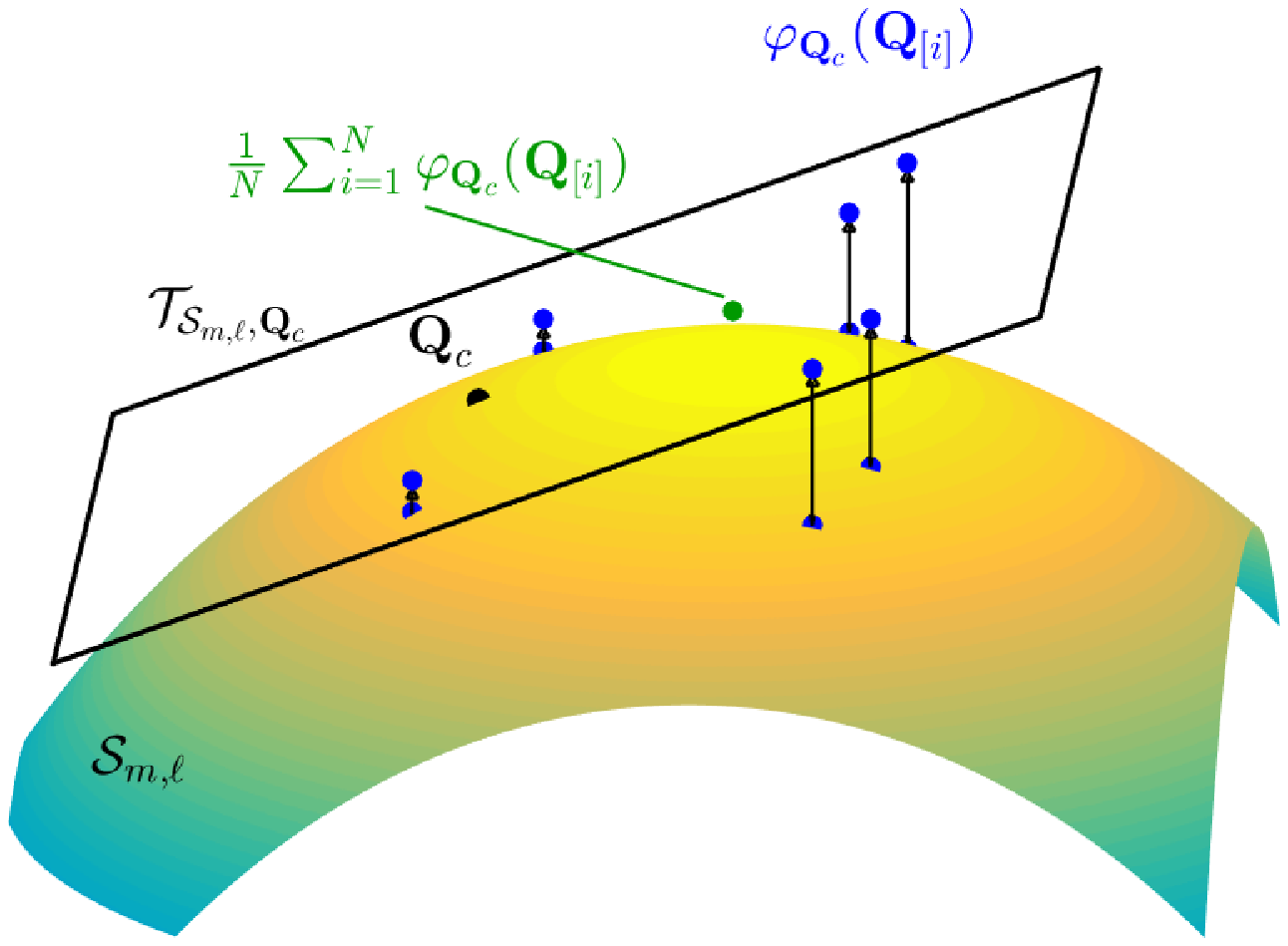}\label{subfig:kn_b}}
  \vspace{-0.4cm}
  \subfloat[ ]{\includegraphics[width=.55\textwidth]{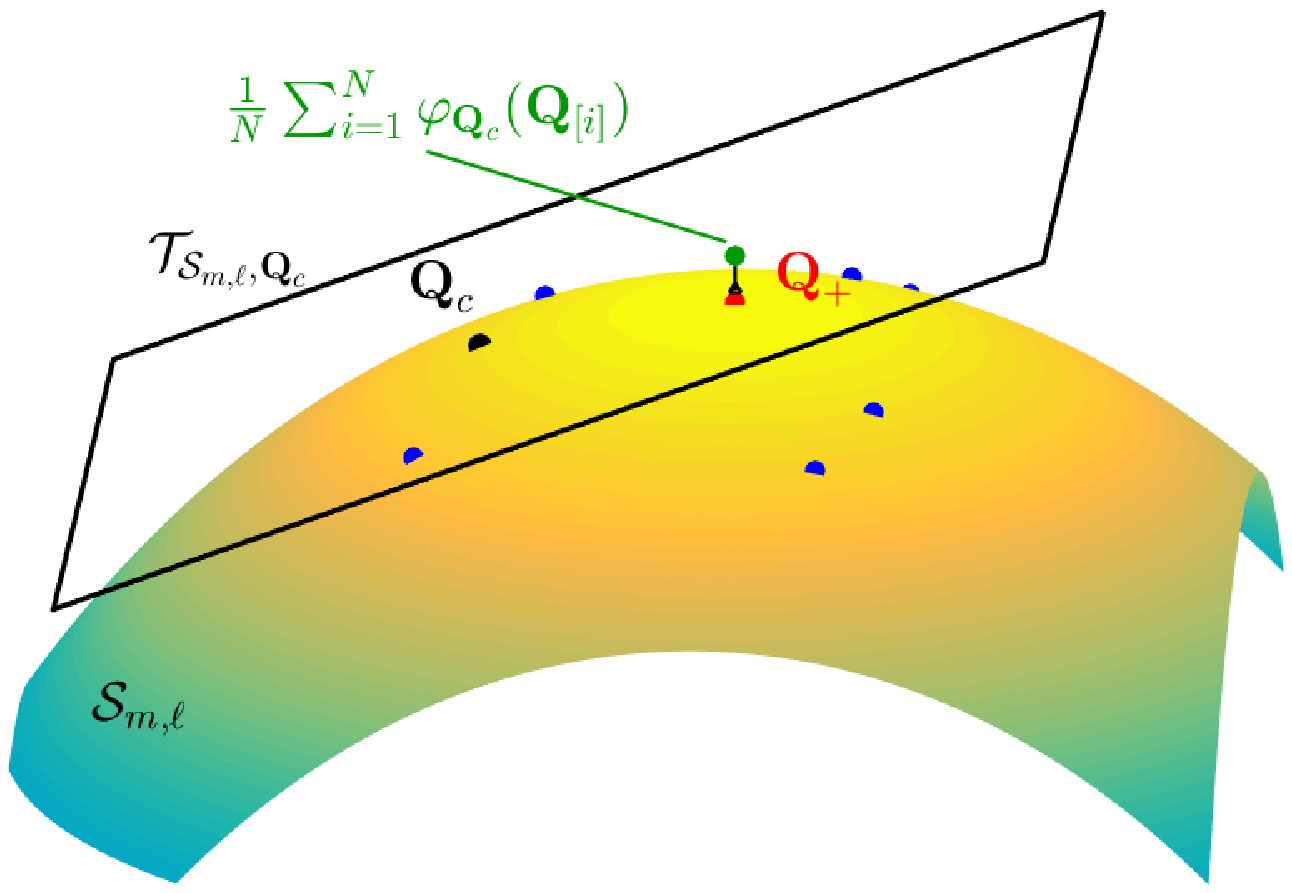}\label{subfig:kn_c}}
  \subfloat[ ]{\includegraphics[width=.55\textwidth]{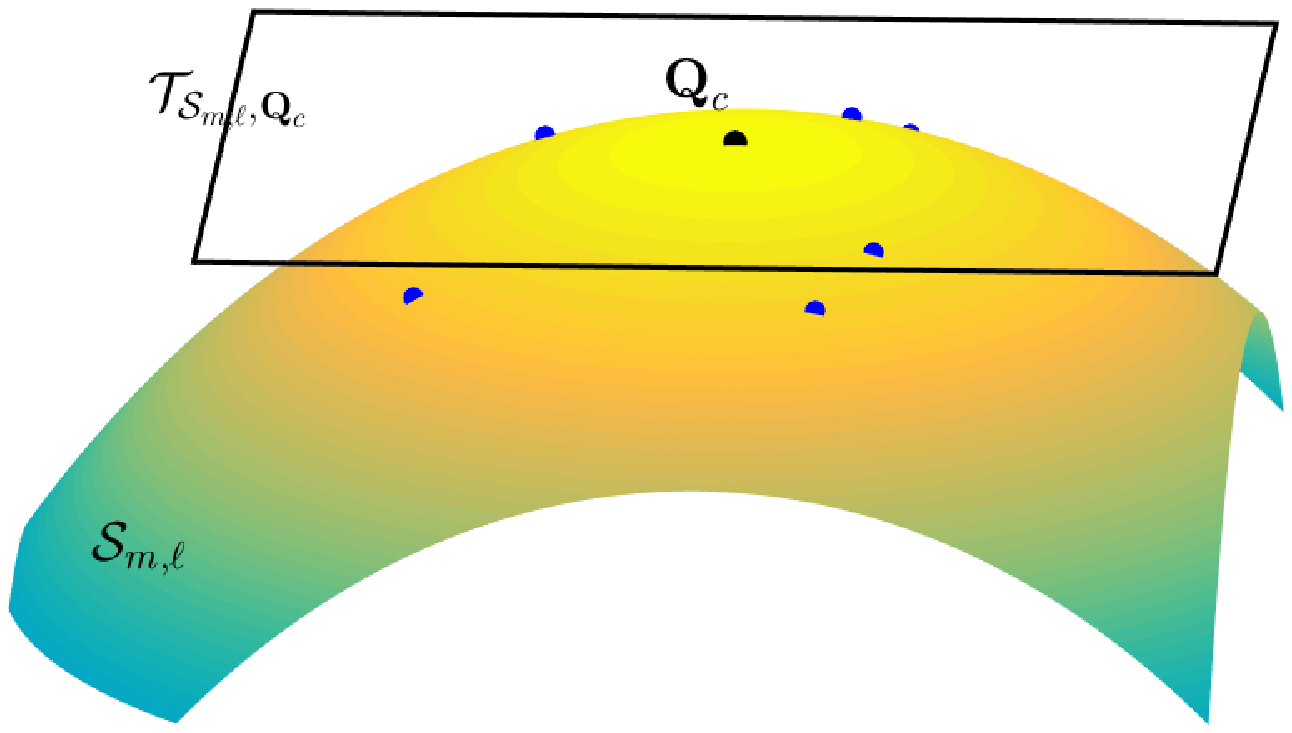}\label{subfig:kn_d}}
  \caption{A conceptual illustration of one iteration of Algorithm~\ref{alg:KN}
    for solving the optimization problem~\eqref{eq:optQ}. (a)~Five $\Qb_\i$ (blue dots)
    on the Stiefel manifold $\St_{m,\ell}$ are lifted to $\varphi_{\Qb_c}(\Qb_\i \Qb_\i^\intercal)$
    on the tangent space ${\mathcal T}_{\Qb_{c}}\St_{m, \ell}$ corresponding to the
    current iterate $\Qb_{c}$ (black dot). (b)~We compute the KN-type average of the
     lifted points, which is identified by the green point.  (c)~The average of the
     lifted points (green dot) is mapped back to the $\Qb_{+}$ (red dot) on the
     Stiefel manifold as the next iterate. (d)~We have moved from $\Qb_{c}$ to
     $\Qb_{+}$ and obtained the new $\Qb_{c}$.}
    \label{fig:kn}
\end{center}
\end{figure}

More details of the algorithm are given below. For the choice of the
initial iterate $\Qb_{\rm ini}$, we select the iterate that has the
largest value of $\tr(\widetilde\Sigmab_\i)$ from the collection
$\{\Qb_\i\}_{i=1}^N$, where $\widetilde\Sigmab_\i$ is the diagonal
matrix consisting of the singular values of $\Yb_\i$ computed in
Step~\ref{line:isvd_proj} of Algorithm~\ref{alg:isvd}. Specifically,
we choose $\Qb_{\rm ini}= \Qb_{[i_{\max}]}$, where
$i_{\max}=\argmax_{i=1,\dots,N} \tr(\widetilde\Sigmab_\i)$. In each
iteration, namely Steps~\ref{step:kn_ave} and~\ref{step:kn_pullback}
of Algorithm~\ref{alg:KN}, we move the current iterate $\Qb_{c}$ to
the next iterate $\Qb_{+}$ via the following procedure.
One iteration of the integration Algorithm~\ref{alg:KN} is
illustrated conceptually in Figure~\ref{fig:kn}.
In particular, Step~\ref{step:kn_ave} of Algorithm~\ref{alg:KN} is
composed of the following two tasks.
\begin{enumerate}
\item As shown in Figure~\ref{fig:kn}(a), we map (or lift) the matrices $\{\Qb_\i\}_{i=1}^N$ to the tangent space ${\mathcal T}_{\Qb_c}\St_{m,\ell}$.
That is, each $\Qb_\i \in \St_{m,\ell}$ is mapped to
\begin{equation}
\varphi_{\Qb_c} (\Qb_\i)=(\Ib_m-\Qb_{c}\Qb_{c}^\intercal) \Qb_\i\Qb_\i^\intercal \Qb_{c} \in {\mathcal T}_{\St_{m, \ell}, \Qb_{c}}.
\label{eq:lifting_Qi}
\end{equation}
Because $\Qb_c^\intercal \{\varphi_{\Qb_c} (\Qb_\i)\} +
\{\varphi_{\Qb_c} (\Qb_\i)\}^\intercal\Qb_c =\0$, we know that
$\varphi_{\Qb_c} (\Qb_\i)$ is indeed a point on ${\mathcal
T}_{\Qb_c}\St_{m,\ell}$ \cite{tagare2011notes}.
\item As shown in Figure~\ref{fig:kn}(b), we then take the average of the mapped matrix points.
Because these mapped matrices are located on ${\mathcal
T}_{\Qb_c}\St_{m,\ell}$, which is a flat space, we can compute
the arithmetic average of $\varphi_{\Qb_c}(\Qb_\i)$
\begin{equation}
\overline{\varphi_{\Qb_c}(\Qb_\i)}
= \frac{1}{N} \sum_{i=1}^{N} \varphi_{\Qb_c} (\Qb_\i)
= (\Ib_m-\Qb_c\Qb_c^\intercal) \overline{\Pb}\Qb_c.
\label{eq:ave_Qi}
\end{equation}
In \eqref{eq:ave_Qi}, we apply \eqref{eq:lifting_Qi} and the
definition of $\overline{\Pb}$ in \eqref{eq:def_P}. This average
is still on the tangent space. Furthermore,
$\overline{\varphi_{\Qb_c}(\Qb_\i)} =\Db_F(\Qb_c)$ by
\eqref{eq:KN_on_tangent}.
\end{enumerate}
In Step~\ref{step:kn_pullback}, as shown in Figure~\ref{fig:kn}(c),
we pull the averaged matrix $\overline{\varphi_{\Qb_c}(\Qb_\i)}$
back to the Stiefel manifold by the inverse map
$\varphi_{\Qb_c}^{-1}$. Specifically,
\begin{equation}
\label{def:KN_Q_next}
\Qb_+
= \varphi_{\Qb_c}^{-1}(\overline{\varphi_{\Qb_c}(\Qb_\i)})
= \Qb_c \Cb+ \overline{\varphi_{\Qb_c}(\Qb_\i)} \Cb^{-1},
\end{equation}
where $\Cb= \left[ \frac{\Ib_\ell}2 + \left[ \frac{\Ib_\ell}4 -
\overline{\varphi_{\Qb_c}(\Qb_\i)}^\intercal
\overline{\varphi_{\Qb_c}(\Qb_\i)} \right]^{1/2}\right]^{1/2}$ by
Theorem~\ref{thm:C} with fixed \mbox{$\tau=1$}.

In short, we move the iterate from $\Qb_{c}$ to $\Qb_{+}$ in the loop of
Algorithm~\ref{alg:KN} by the following procedure. (i) $\Qb_\i$ are
mapped to ${\mathcal T}_{\Qb_c}\St_{m, \ell}$ by
$\varphi_{\Qb_{c}}$, (ii) the mapped matrices are averaged as
$\overline{\varphi_{\Qb_c}(\Qb_\i)} = \Db_F(\Qb_c)$, and finally,
(iii) $\Db_F(\Qb_c)$ is mapped back to the manifold by the inverse
map $\varphi_{\Qb_c}^{-1}$ to obtain the next iterate $\Qb_+$. This
process can be summarized in one line: $\Qb_+\leftarrow
\varphi_{\Qb_c}^{-1} \left(\frac 1{N}\sum_{i=1}^N
\varphi_{\Qb_c}(\Qb_\i)\right)$.

\subsection{Convergence and Remarks}
\label{sec:isvd_conv}

Algorithm~\ref{alg:KN} is a fixed-point iteration with step size
$\tau=1$. The update from the current $\Qb_c$ to the next $\Qb_+$
can be written as
\begin{equation}\label{fixed_point}
\Qb_+ = g(\Qb_c) = \Qb_c\Cb + \Xb\Cb^{-1},
\end{equation}
where $\Cb$ and $\Xb$ both depend on $\Qb_c$ and can be denoted as
$\Cb(\Qb_c)$ and $\Xb(\Qb_c)$, respectively. Recall that
Algorithm~\ref{alg:KN} is used to find the maximizer of the
objective function $F(\Qb) = \frac12
\tr\left(\Qb^\intercal\overline{\Pb}\Qb\right)$. Let $\Qb_*$ consist
of the leading $\ell$ eigenvectors of $\overline{\Pb}$.
Specifically, $\Qb_*$ consists of the maximizer (uniquely up to an
orthogonal transformation) of the objective function $F(\Qb)$.
Further, let ${\cal N}_\varepsilon(\Qb_*) := \{\Qb\in\St_{m,\ell}:
\|\Qb- \Qb_*\|_F<\varepsilon\}$ be an $\varepsilon$-neighborhood of
$\Qb_*$ in $\St_{m,\ell}$. We can see that $\Qb_*$ is a fixed point
for $g$ in~(\ref{fixed_point}). We establish the convergence for the
fixed-point iteration in Theorem~\ref{thm:fixed_point}. The theorem
suggests that Algorithm~\ref{alg:KN} converges if it starts from an
initial iterate that belongs the $\varepsilon$-neighborhood of an
equivalent version of $\Qb_*$. The equivalence is in the sense of an
orthogonal transformation multiplied on the right side of $\Qb_*$.
The proof of Theorem~\ref{thm:fixed_point} is given in
Appendix~\ref{a.thm.fixed_point}.
\begin{thm}\label{thm:fixed_point}
There exists an $\varepsilon>0$ such that Algorithm~\ref{alg:KN}
converges, provided that the iteration starts from an initial
$\Qb_{\rm ini}\in{\cal N}_\varepsilon(\Qb_* \Rb_0)$, where $\Rb_0$ is
an arbitrary orthogonal matrix.
\end{thm}


We conclude the discussion of Algorithm~\ref{alg:KN} with the
following remarks. First, we bridge the theoretical aspect of the
optimal representation $\overline{\Qb}$ and the numerical scheme
shown in Algorithm~\ref{alg:KN}. Because $\overline{\Qb}$ is the
solution of the optimization problem~\eqref{eq:optQ}, we have
$\Db_F(\overline{\Qb})=\0$. Equation~\eqref{eq:KN_on_tangent}
further suggests that $\frac 1N\sum_{i=1}^N
\varphi_{\overline{\Qb}}(\Qb_\i )=\0$. Moreover, by the definition
\eqref{KN_ave}, we can obtain the Kolmogorov-Nagumo-type average of
$\Qb_\i$ in terms of $\overline{\Qb}$:
\begin{equation}
\varphi_{\overline{\Qb}}^{-1} \left(\frac 1{N}\sum_{i=1}^N\varphi_{\overline{\Qb}}(\Qb_\i)\right)
= \varphi_{\overline{\Qb}}^{-1} \left(\0 \right)
= \overline{\Qb}.
\label{eq:ave_Qbar}
\end{equation}
The last equality holds because of the following. For $\Xb=\0$,
Equation~\eqref{eq:matrixC} suggests that
\begin{equation}
\Cb=\Ib_\ell,
\label{eq:CeqI}
\end{equation} and Equation~\eqref{eq:phi_inv}
further suggests that $\varphi_{\overline{\Qb}}^{-1} \left(\0
\right) = \overline{\Qb}$. Equation~\eqref{eq:ave_Qbar} indicates
that the optimal representation $\overline{\Qb}$ is a
Kolmogorov-Nagumo type average and also a fixed point in
Algorithm~\ref{alg:KN} with the corresponding projected gradient
equal to zero. These facts represent a theoretical background for
computing $\overline{\Qb}$,  and Algorithm~\ref{alg:KN} provides a
numerical method to compute $\overline{\Qb}$.

Second, we use small $\| \Cb - \Ib_\ell\|$ as the stopping criterion
on Step~\ref{ln:KN_conv} of Algorithm~\ref{alg:KN} based on the fact
shown in \eqref{eq:CeqI}. This choice of stopping criterion can be
viewed from the small change between $\Qb_+$ and $\Qb_c$. When $\Cb$
is close to the identity matrix, $\Qb_+$ is close to $\Qb_c$. It is
worth mentioning that $\Cb$ is a small matrix with dimensions $\ell
\times \ell$ and is computed in the iteration of
Algorithm~\ref{alg:KN}. Therefore, the stopping criterion does not
require extra computational effort.

Third, the constrained maximization problem \eqref{eq:optQ} can be
solved using the gradient ascent method proposed in
\cite{wen2013feasible}. The method starts from an initial $\Qb_{\rm
ini}\in \St_{m, \ell}$ and updates the current iterate $\Qb_{c}$ by
searching the next iterate $\Qb_{+}$ on a curve lying on the Stiefel
manifold $\St_{m, \ell}$ to satisfy the orthogonality constraint.
The curve is obtained by mapping the projected gradient defined in
\eqref{def:DF01} to the Stiefel manifold via a Cayley transform. An
efficient step size selection along the curve can accelerate the
overall convergence. On the other hand,   Theorem~\ref{thm:C}
presents a curve along the direction of the projected gradient on
the manifold. In Algorithm~\ref{alg:KN}, it is equivalent to setting
the step size as $\tau=1$, and Algorithm~\ref{alg:KN} can
consequently be viewed as a gradient ascent method.

We have proposed and analyzed Algorithm~\ref{alg:KN} to compute
$\overline{\Qb}$ by solving the constrained optimization
\eqref{eq:optQ} (and \eqref{eq:QiQ_1} equivalently). With the
computed $\overline{\Qb}$, we can use iSVD, i.e.,
Algorithm~\ref{alg:isvd}, to perform the approximate SVD defined in
\eqref{eq:isvd} with multiple random sketches. In the next section,
the iSVD is analyzed statistically.

\section{Statistical Analysis}
\label{sec:stat}

In this section, we present some theoretic statistical analysis on
iSVD. First, we prove a Strong Law of Large Numbers (SLNN) result in
Theorem~\ref{thm:lln} to show that iSVD \eqref{eq:isvd} can perform
as well as the full data SVD \eqref{eq:rank_k_svd} as the number of
random sketches $N$ goes to infinity. Next, consistencies in terms
of subspace and singular vectors are asserted in
Theorem~\ref{thm:cnstcy}. Finally, we determine a rate of
convergence by the Central Limit Theorem (CLT) in
Theorem~\ref{thm:clt}.


{\bf Strong Law of Large Numbers.} From Theorem~\ref{lemma:key} and
the fact that the absolute values of entries of a projection matrix
are bounded by one, we have the following immediate result based on
Theorem~\ref{lemma:key}.
\begin{thm}[Strong Law of Large Numbers]\label{thm:lln}
We have
\[\lim_{N\to\infty} \frac1N\sum_{i=1}^N \Qb_{[i]}\Qb_{[i]}^\intercal
= \Ub\Lambdab\Ub^\intercal ~~\mbox{with probability one},\] where
$\Lambdab$ is given in Theorem~\ref{lemma:key} and $\Ub$ is the true
left singular vectors of the underlying matrix $\Ab$ in decreasing
order.
\end{thm}


{\bf Consistency.} Next, we establish the consistency between the
left singular vectors computed by iSVD and the true left singular
vectors. We prove Lemma~\ref{lemma:approx_ell} first. Based on the
lemma, we prove the consistency in Theorem~\ref{thm:cnstcy}. See Appendix~\ref{proof.consistency} for the proofs of the lemma and the theorem.

\begin{lem}\label{lemma:approx_ell}
Let $\Ub=\left[\ub_1,\ldots,\ub_m\right]$ be an arbitrary point in
$\St_{m,m}$, and let $\Lambdab$ be a diagonal matrix with decreasing
diagonal entries $1> \lambda_1 \geq \lambda_2 \geq \ldots
\lambda_\ell > \lambda_{\ell+1}\ge\ldots\ge \lambda_m \geq0$.
Consider the following minimization problem:
\[\Qb_\opt =\argmin_{\Qb\in \St_{m,\ell}} \left\| \Ub\Lambdab\Ub^\intercal
 - \Qb \Qb^\intercal\right\|_F^2.\]
Then, we have $\Qb_\opt \Qb_\opt^\intercal =\Ub_\ell
\Ub_\ell^\intercal$, where $\Ub_\ell
=\left[\ub_1,\ldots,\ub_{\ell}\right]$.
\end{lem}

\begin{thm}[Consistency of subspaces and singular vectors.]\label{thm:cnstcy}
Assume that the diagonal entries of $\Sigmab$ (i.e., singular values
of $\Ab$) satisfy the condition:
$\sigma_1>\sigma_2>\cdots>\sigma_{\ell}>\sigma_{\ell+1}
\ge\cdots\ge\sigma_{m}\ge0$. Then, we have the following properties.
(a)~$\lim_{N\to \infty}
\overline{\Qb}\,\overline{\Qb}^\intercal=\Ub_{\ell}
\Ub_{\ell}^\intercal$ with probability one. (b)~Let
$\widehat{\Wb}_{\ell}$ consist of the left singular vectors of
$\overline{\Qb}^\intercal\Ab$, and let
$\widehat{\Ub}_{\ell}=\overline{\Qb}\widehat{\Wb}_{\ell}$ as
described in Algorithm~\ref{alg:isvd}. Then, for any $j\leq \ell$,
we have
\[
\lim_{N\to \infty} \left|\widehat{\ub}_j^\intercal  \ub_j\right|=1 \quad
\mbox{with probability one},
\]
where $\widehat{\ub}_j$  is the $j${\rm th} column of
$\widehat{\Ub}_{\ell}$ and ${\ub}_j$  is the $j${\rm th} column of
${\Ub}_\ell$.
\end{thm}

Note that the consistency established in Theorem~\ref{thm:cnstcy} is
valid for the entire $\widehat\Ub_\ell$, where $\ell$ is the
sampling dimension. However, we expect a more accurate
$\widehat\Ub_k$ using a larger sampling dimension $\ell$. See
Table~\ref{tab:notation} for the definitions of $k$ and $\ell$.

{\bf Central Limit Theorems.} Because $\Qb_\ione$, $\Qb_\itwo$,
$\ldots$, $\Qb_\iN$ are i.i.d., so are
$\Qb_\ione\Qb_\ione^\intercal$, $\Qb_\itwo\Qb_\itwo^\intercal$,
$\ldots$, $\Qb_\iN\Qb_\iN^\intercal$; and they have finite second
moments. The following theorem is an immediate CLT result from
Theorem~\ref{lemma:key}.
\begin{thm}[Central Limit Theorem I]\label{thm:clt}
We have
\begin{equation}\label{eq:clt}
\frac1{\sqrt N}\sum_{i=1}^N \vec\left(\Qb_{[i]}\Qb_{[i]}^\intercal
- \Ub\Lambdab\Ub^\intercal\right) \stackrel{d}{\rightsquigarrow}
  {\cal N}(\0,\Tb_1),~~ {\rm as~} N\to\infty,
\end{equation}
where $\Tb_1$ is a certain positive definite matrix.
\end{thm}

Theorem~\ref{thm:clt} is a CLT on the average of projection
matrices. However, a more sensible CLT should be for the
singular vectors $\widehat \ub_j$ estimated by iSVD.
Note that $\widehat \ub_j$ (or $\ub_j$) is a function of
$\frac1N\sum_{i=1}^N\Qb_{[i]}\Qb_{[i]}^\intercal$ (or
$\Ub\Lambdab\Ub^\intercal$). By the delta-method
to~(\ref{eq:clt}), we can establish the following CLT on $\widehat
\ub_j$. See Appendix~\ref{proof.clt} for the proof.

\begin{thm}[Central Limit Theorem I\!I]\label{thm:clt2}
We have
\[ {\sqrt N} \left(\widehat\ub_j - \ub_j\right)\stackrel{d}{\rightsquigarrow}
 {\cal N}(\0,\Tb_2),~~ {\rm as~} N\to\infty,\] where $\Tb_2 =
\Deltab_j\Tb_1 \Deltab_j^\intercal$ and $\Deltab_j =
\frac{\partial\ub_j}{\partial\vec
(\Ub\Lambdab\Ub^\intercal)^\intercal}$ is given by~(\ref{Delta})
below.
\end{thm}

From Theorem~\ref{thm:clt2}, we know that $\widehat\ub_j- \ub_j
=O_p(N^{-1/2})$ and so is $\widehat\Ub_\ell- \Ub_\ell
=O_p(N^{-1/2})$. Then,
$\|\widehat\Ub_\ell\,\widehat\Ub_\ell^\intercal - \Ub_\ell
\Ub_\ell^\intercal\|_F^2=O_p(N^{-1})$ and
$\|\widehat\Ub_\ell\,\widehat\Ub_\ell^\intercal - \Ub_\ell
\Ub_\ell^\intercal\|_{\rm sp}=O_p(N^{-1/2})$. Note that
$\widehat\Ub_\ell\,\widehat\Ub_\ell^\intercal  =
\overline{\Qb}\,\overline{\Qb}^\intercal$. From $\|\Ub_\ell
\Sigmab_\ell  \Vb_\ell -\Ab \|_F^2 =\sigma_{k+1}^2+\ldots
+\sigma_m^2$ and $\|\Ub_\ell \Sigmab_\ell \Vb_\ell  -\Ab \|_{\rm
sp}=\sigma_{k+1}$, we have
\begin{equation}\label{ourbound_F}
\|\overline{\Qb}\,\overline{\Qb}^\intercal \Ab -\Ab \|_F^2
=\sigma_{k+1}^2+\ldots \sigma_m^2+O_p(N^{-1})
\end{equation}
and
\begin{equation}\label{ourbound}
\|\overline{\Qb}\,\overline{\Qb}^\intercal \Ab -\Ab \|_{\rm sp}^2
 =\sigma_{k+1}^2+O_p(N^{-1}).
 \end{equation}
Specifically, as $N\to\infty$, we can achieve tight bounds in both
the Frobenius norm and the spectral norm by integrating multiple
random sketches.

\section{Numerical Results}
\label{sec:num}

We conduct numerical experiments to study the performance of the
proposed algorithms.
To test the proposed iSVD, we construct the following test matrices,
which are similar to the test matrices used in
\cite{rokhlin2009randomized}. Let the matrix
$\Ab= \Hb_d\, \Sigmab\, \Hb_{d+1}^\intercal \in
\mathbb{R}^{2^d\times 2^{d+1}},$
where $\Hb_d$ is the Hadamard matrix of size $2^d\times 2^d$ and
$\Sigmab$ is a diagonal matrix of size $2^d\times 2^{d+1}$. Note
that, for a Hadamard matrix, $\Hb_{d}^\intercal = \Hb_d$ and
$\Hb_{d}^\intercal \Hb_d =\Ib_{2^d}$. Let the desired rank be
$k=10$. We set the $j$th diagonal entry of $\Sigmab$ as follows:
\begin{equation}\label{eq:Sigma}
\Sigmab_{j,j} = \sigma_{j}
= \left\{\begin{array}{ll}
\sigma_1^{\lfloor j/2\rfloor / 5}, & j=1,3,5,7,9,\\[1.1ex]
1.5\sigma_{j+1} & j=2,4,6,8, 10,\\[0.8ex]
0.001, & j=11,\\[0.8ex]
\sigma_{11}\cdot\frac{m-j}{m-11}, & j=12,\dots,m.
\end{array}\right.
\end{equation}
Here, $\lfloor j/2 \rfloor$ is the greatest integer less than or
equal to $j/2$. Our $\Sigmab$ is modified
from~\cite{rokhlin2009randomized} to distinguish the singular
values, and thus, individual singular vectors can be uniquely
identified. Note that $\Hb_d$ is an orthogonal matrix. Thus, the SVD
of the test matrix $\Ab$ is known to be $\Hb_d\, \Sigmab\,
\Hb_{d+1}^\intercal$, where the columns of $\Hb_d$ and $\Hb_{d+1}$
are the left and right singular vectors, respectively, and
$\sigma_j$ are singular values.

The experimental settings are $d=9,11,13,15,17,19$, $k=10$, $p=12$,
$\ell=k+p=22$, $q=0,1$, and $N=10,50,100, 200$. For an initial
$\Qb_{\rm ini}$, we select from the collection $\{\Qb_\i\}_{i=1}^N$.
We choose the $\Qb_\i$ that has the largest value of
$\tr(\widetilde\Sigmab_\i)$, where $\widetilde\Sigmab_\i$ is the
diagonal matrix consisting of the singular values of $\Yb_\i$
computed in Step~\ref{line:isvd_proj} of Algorithm~\ref{alg:isvd}.
To evaluate the accuracy of approximate SVD, we use the following
similarity for comparing the computed and true leading $k$
individual singular vectors. Recall that $\widetilde \Ub_k = \{
\widetilde{\ub}_1, \ldots, \widetilde{\ub}_j, \ldots,
\widetilde{\ub}_k \}$ and $\widehat \Ub_k=\{\widehat{\ub}_1, \ldots,
\widehat{\ub}_j, \ldots, \widehat{\ub}_k\}$ consist of the rank-$k$
left singular vectors computed by Algorithm~\ref{alg:rsvd} (rSVD)
and Algorithm~\ref{alg:isvd} (iSVD), respectively. $\Ub_k = \{
\ub_1, \ldots, \ub_j, \ldots, \ub_k \}$ consists of the true left
singular vectors of $\Ab$. We measure the similarity between the
$j$th computed singular vector $\widetilde \ub_j$ (or $\widehat
\ub_j$) and the true singular vector $\ub_j$ by computing
$|\widetilde \ub_j^\intercal \ub_j|$ (or $|\widehat \ub_j^\intercal
\ub_j|$) for $j=1,\ldots,k$. If the computed singular vector has no
error, then $|\widetilde \ub_j^\intercal \ub_j|=1$ (or $|\widehat
\ub_j^\intercal \ub_j|=1$). Note that we present only the results
regarding the left singular vectors. The results involving the right
singular vectors are similar and ignored here.
Algorithm~\ref{alg:KN} is stopped if $\|\Cb-\Ib_\ell\|_F$ is less
than $10^{-5}$. The numerical experiments are conducted on a
workstation equipped with an Intel E5-2650 v3 CPU (with a 25 MB
cache and 2.30 GHz clock rate) and 256 GB of main memory. The
algorithms are implemented in MATLAB version 2015b.

\begin{figure}
\begin{center}
\subfloat[$q=0$]{\includegraphics[width=.55\textwidth]{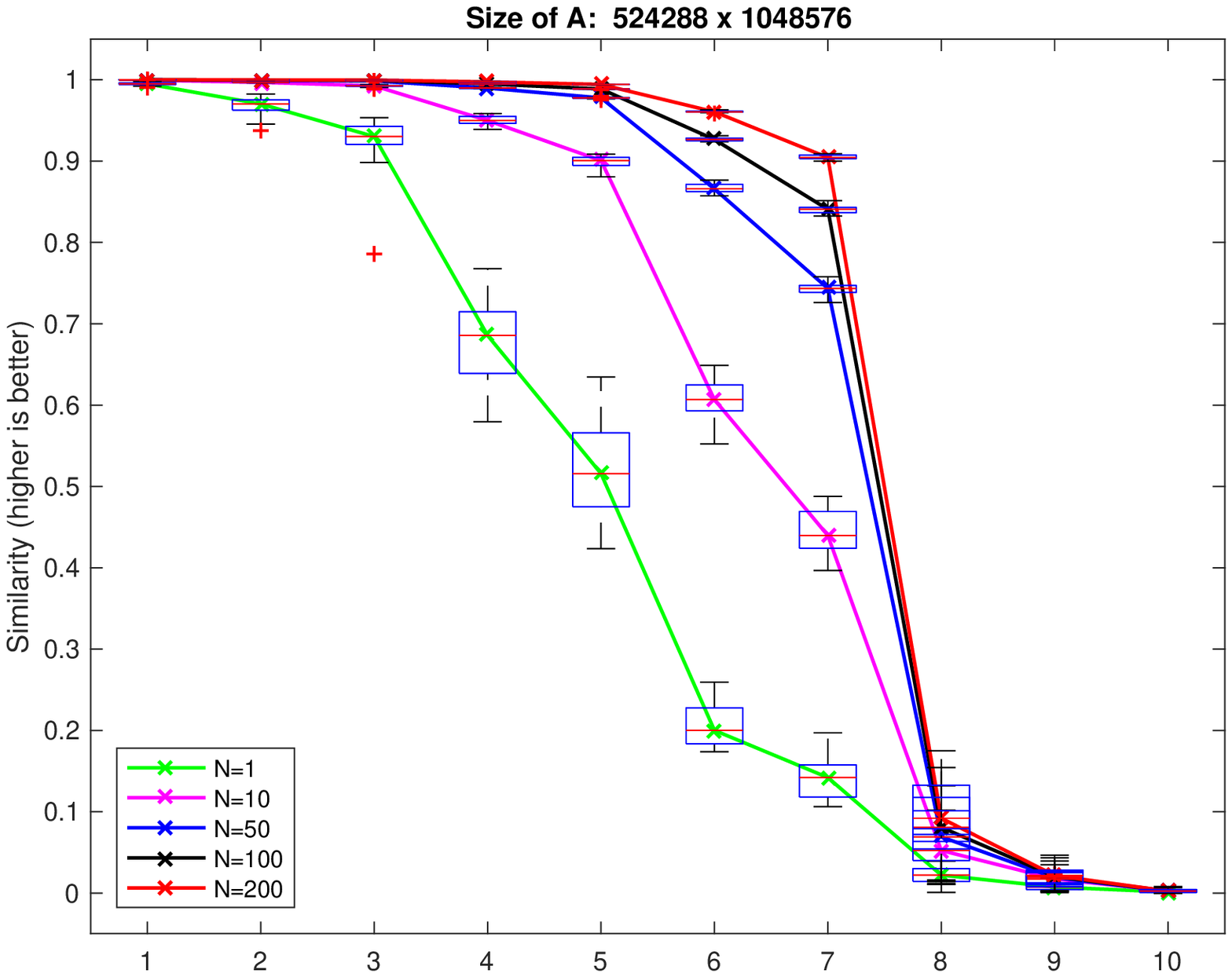}\label{subfig:h19_q0}}
\subfloat[$q=1$]{\includegraphics[width=.595\textwidth]{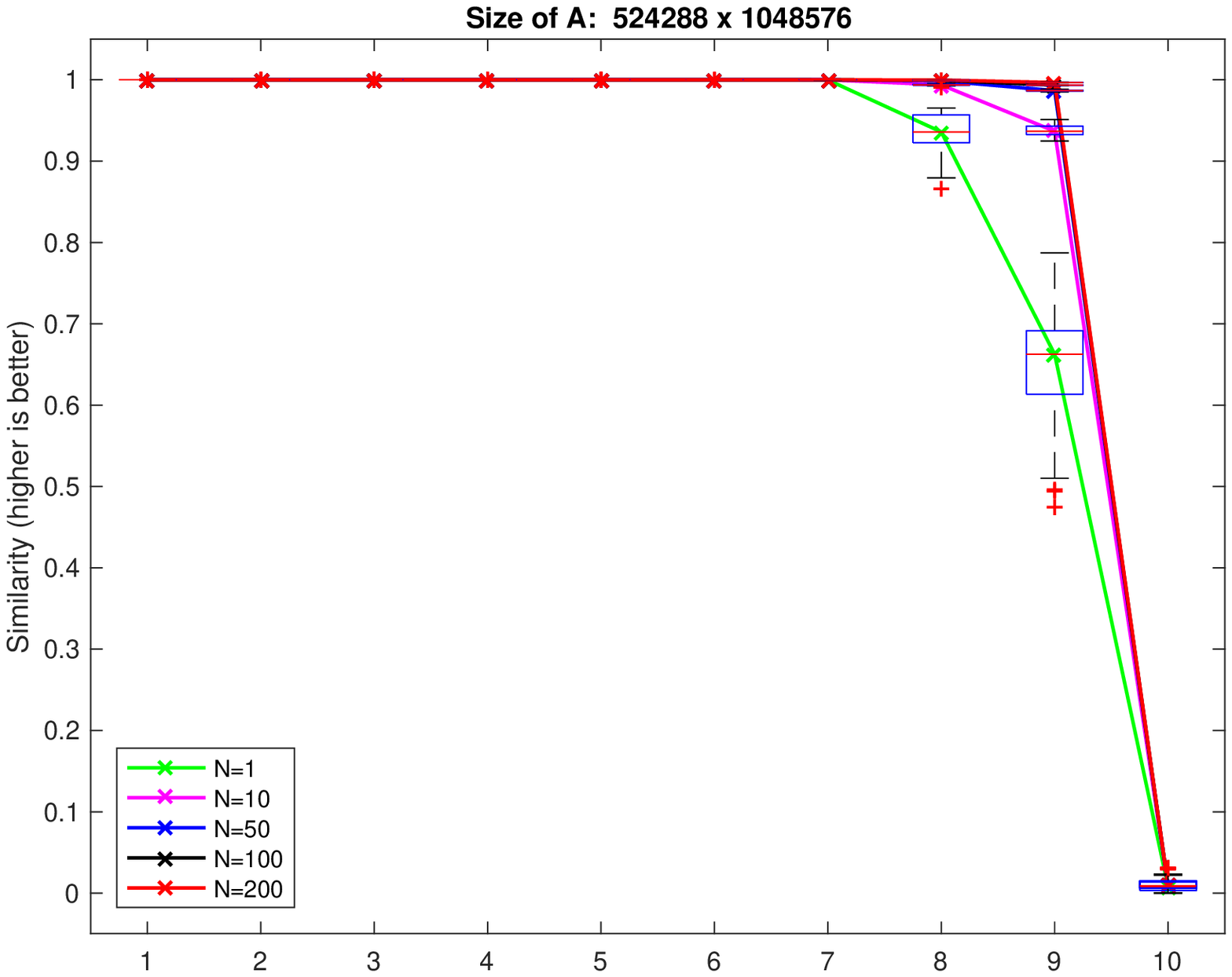}\label{subfig:h19_q1}}
\end{center}
\caption{The similarity results for $d=19$, $\ell=22$, $q=0,1$, and $N=$1, 10, 50, 100,  200.
The figure also shows the box plots indicating the median and the 25th
and 75th percentiles of the similarities out of $30$ replicated runs.
Higher similarities (up to $1$) are better.}
\label{fig:similarity}
\end{figure}

We report the accuracies and variations in the computed singular
values and singular vectors in Figure~\ref{fig:similarity} and
Table~\ref{tab:dallq0q1} using different parameters. We highlight
the following observations.
\begin{itemize}
\item
{\it The similarity (accuracy) of the singular vectors increases
as the number of random sketches $N$ increases.} For each
singular vector, we examine the accuracy performance in terms of
the similarity between the computed and true singular vector.
Figure~\ref{fig:similarity} shows the singular vector similarity
results with box plots. In the figure, the matrix size is
$2^{19}\times 2^{20}$ ($d=19$), the sampling dimension
$\ell=22$, the number of random sketches $N=$1, 10, 50, 100, and
200, and the exponent of the power method in
Step~\ref{step:randproj} of Algorithms~\ref{alg:rsvd}
and~\ref{alg:isvd} (i.e., $q$) equals $0$ or $1$. Higher
similarities (up to $1$) are better. It is clear that larger $N$
results in higher similarity in all  the tested cases. Some of
the improvements can be significant, especially for several
cases when $q=0$ and the $9$th singular vector for $q=1$. Note
that the $10$th singular vector is difficult to compute. This is
because the $10$th eigenvalue belongs to a cluster of singular
values, and it is difficult to distinguish the singular vectors
of these slow-decaying singular values.

\item {\it The rank-$k$ matrix error decreases as the number of
random sketches $N$ increases.} We also examine the accuracy
performance with respect to the combination of the singular
values and singular vectors. In particular, we compute the
rank-$k$ matrix error $\| \Ub_k\Sigmab_k \Vb^\intercal_k-
\widehat{\Ub}_k \widehat{\Sigmab}_k \widehat{\Vb}^\intercal
_k\|_{F}$. This error evaluates the difference between the
estimated and true rank-$k$ SVD. We experiment with different
$d$ to better present the trend of the integration effect for
 $d$. The observations hold for all the experiments for $d=$9,
11, 13, 15, 17, and 19 with $q=0$ and $q=1$, as shown in
Table~\ref{tab:dallq0q1}.

\item
{\it Overall, the stochastic variation in similarity of a
singular vector to its target decreases as the number of random
sketches $N$ increases.} This welcomed result can be expected
because more random sketches have been integrated, and thus, the
averaged sketch becomes more stable and with less stochastic
variation. Such an observation holds for almost all the
numerical results shown in Table~\ref{tab:dallq0q1}.
\end{itemize}

Furthermore, we investigate the effect of increasing the sampling
dimension $\ell$ for rSVD ($N=1$) and compare the results with iSVD
($N>1$), which uses $\ell=22$ and $N=200$, resulting in $\ell\times
N = 4400$ samples in total. In these numerical experiments, $d=19$,
\mbox{$q=0$}, and the number of replicated runs is $30$.
Table~\ref{tab:comp} shows the Frobenius norm of the error matrix
(i.e., $\| \Ub_k\Sigmab_k \Vb_k^\intercal - \widetilde{\Ub}_k
\widetilde{\Sigmab}_{k}\,\widetilde{\Vb}_{k}^\intercal \|_{F}$ for
rSVD and $\| \Ub_k\Sigmab_k \Vb_k^\intercal - \widehat{\Ub}_k
\widehat{\Sigmab}_k \widehat{\Vb}^\intercal _k\|_{F}$ for iSVD). Two
main observations are highlighted below.
\begin{itemize}
\item
{\it In rSVD ($N=1$), a larger $\ell$ results in smaller average
errors and smaller standard deviations.} This observation is
reasonable because when we sketch a  greater number of sampling
dimensions, more information of the leading singular vectors is
collected.
\item {\it SVD computed by iSVD with smaller sampling dimensions (via multiple sketches) is better than rSVD with large sampling dimensions (via a single sketch).} We
compare the result obtained by iSVD with $\ell=22$ and $N=200$
($4400$ sampling dimensions in total) with the results obtained
by rSVD with various $\ell$ and $N=1$. As shown in
Table~\ref{tab:comp}, iSVD outperforms rSVD in all cases except
for the case with $\ell=4400$. For the case in which
$\ell=4400$, rSVD performs slightly better. This observation
suggests the advantage of integration. In addition, even without
adopting parallelism, taking $200$ random sketches with
$\ell=22$ and integrating them is relatively efficient compared
to executing an rSVD with $\ell=3000$ in terms of both precision
and time.

\end{itemize}

\begin{table}
\centering
\footnotesize 
\begin{tabular}{cccccccccc|c|c|c|c|c|}\hline
d &N=1&N=10&N=50&N=100&N=200\\\hline
\multicolumn{6}{c}{(a) $q=0$}\\ \hline
9&    1.04e-02    (6.56e-04)&    3.79e-03    (1.18e-04)&    1.74e-03    (4.37e-05)&    1.23e-03    (2.46e-05)&    8.71e-04    (1.52e-05)\\
11&    1.89e-02    (1.12e-03)&    6.74e-03    (1.51e-04)&    3.25e-03    (5.94e-05)&    2.32e-03    (3.07e-05)&    1.67e-03    (1.90e-05)\\
13&    3.49e-02    (2.69e-03)&    1.22e-02    (2.44e-04)&    5.83e-03    (5.45e-05)&    4.32e-03    (2.69e-05)&    3.30e-03    (1.63e-05)\\\hline
15&    6.20e-02    (4.20e-03)&    2.21e-02    (4.15e-04)&    1.06e-02    (9.62e-05)&    7.78e-03    (3.94e-05)&    5.72e-03    (2.11e-05)\\
17&    1.12e-01    (5.76e-03)&    4.03e-02    (8.24e-04)&    1.95e-02    (1.72e-04)&    1.44e-02    (7.94e-05)&    1.09e-02    (5.12e-05)\\
19&    1.92e-01    (1.26e-02)&    7.14e-02    (1.59e-03)&    3.52e-02    (2.93e-04)&    2.60e-02    (1.61e-04)&    1.95e-02    (6.30e-05)\\\hline
\multicolumn{6}{c}{(b) $q=1$}\\ \hline
9&    1.08e-03    (1.50e-04)&    4.30e-04    (3.55e-05)&    1.95e-04    (1.40e-05)&    1.37e-04    (8.63e-06)&    9.75e-05    (5.96e-06)\\
11&    1.53e-03    (1.03e-04)&    7.61e-04    (3.89e-05)&    3.68e-04    (1.43e-05)&    2.62e-04    (1.06e-05)&    1.87e-04    (6.84e-06)\\
13&    1.83e-03    (5.54e-05)&    1.23e-03    (4.43e-05)& 6.89e-04
(2.00e-05)&    5.05e-04    (1.32e-05)&    3.65e-04
(7.55e-06)\\\hline
15&    2.14e-03    (9.61e-05)&    1.64e-03    (4.45e-05)&    1.17e-03    (2.59e-05)&    9.27e-04    (1.61e-05)&    7.22e-04    (1.26e-05)\\
17&    2.97e-03    (2.81e-04)&    1.93e-03    (2.30e-05)&    1.78e-03    (7.72e-06)&    1.75e-03    (8.45e-06)&    1.74e-03    (9.44e-06)\\
19&    4.14e-03    (2.77e-04)&    2.35e-03    (5.20e-05)& 1.89e-03
(6.91e-06)&    1.82e-03    (2.81e-06)&    1.78e-03
(9.79e-07)\\\hline
\end{tabular}

\caption{The average norm and standard deviation (in parentheses) of
the error matrices (i.e.,  $\| \Ub\Sigmab \Vb^\intercal-
\widehat{\Ub}_k \widehat{\Sigmab}_k \widehat{\Vb}^\intercal
_k\|_{F}$) over 30 replicated runs with various matrices of size
$2^{d}\times 2^{d+1}$. The sampling dimension $\ell=22$ and the
exponent $q=0$ or $q=1$.} \label{tab:dallq0q1}
\end{table}

\begin{table}
\centering
\footnotesize 
\begin{tabular}{cccccccccc}\cline{1-3}\cline{5-7}
$\ell$ & $N$ (Alg.) & Ave (std) of errors &\ &
$\ell$ & $N$ (Alg.) & Ave (std) of errors \\ \cline{1-3}\cline{5-7}
  22 & 1 (rSVD) & 1.90e-01 (1.39e-02) &  & 3000 &   1 (rSVD) & 2.06e-02 (9.16e-05) \\
 500 & 1 (rSVD) & 4.63e-02 (6.27e-04) &  & 4400 &   1 (rSVD) & 1.73e-02 (6.27e-05) \\
1000 & 1 (rSVD) & 3.40e-02 (2.84e-04) &  &   22 & 200 (iSVD) & 1.95e-02 (6.30e-05) \\ \cline{1-3}\cline{5-7}
\end{tabular}
\caption{We use rSVD (with various sampling dimensions $\ell$ and
$N=1$) and iSVD (with $\ell=22$ and $N=200$) to compute the first
$10$ singular values and singular vectors. The table shows the
averages (and standard deviations) of the error matrix norms (i.e.,
$\| \Ub\Sigmab \Vb^\intercal - \widetilde{\Ub}_k
\widetilde{\Sigmab}_{k}\,\widetilde{\Vb}_{k}^\intercal \|_F$ for
rSVD and $\| \Ub\Sigmab \Vb^\intercal - \widehat{\Ub}_k
\widehat{\Sigmab}_k \widehat{\Vb}^\intercal _k\|_{F}$ for iSVD) out
of $30$ replicated runs. The matrix size is $2^{19}\times 2^{20}$,
and the power exponent $q=0$.}\label{tab:comp}
\end{table}

\section{Conclusions}
\label{sec:con}

We have proposed and analyzed a Monte Carlo-type algorithm for
computing the rank-$k$ SVD of large matrices. The proposed algorithm
integrates multiple leading low-dimensional subspaces projected by
multiple random sketches. The integrated subspace is the solution of
the optimization problem constrained by the matrix Stiefel manifold
that best represents the multiple random projected subspaces. To
solve the optimization problem, we propose an iterative method based
on the Kolmogorov-Nagumo-type average of the multiple subspaces.
Theoretical analyses reveal the insights of the proposed algorithms.
Numerical experiments suggest that the integrated SVD can achieve
higher accuracy and less stochastic variation in singular vectors
using multiple random sketches.

It is interesting to generalize iSVD to other problems. First, we
plan to investigate how  iSVD performs if we replace the Gaussian
random projections by the column random sampling. Unlike the
Gaussian random projections, which involve matrix-matrix
multiplications $\Ab^\intercal \Ab \Omegab$, the random column
sampling can be implemented by column extractions without involving
matrix-matrix multiplications, and therefore leads to a significant
savings in computational time and memory usage, especially for
large-scale matrices. However, the sketched subspaces contain less
information about the leading subspaces, which may decrease the
accuracy, and some statistical properties are different from the
cases in Gaussian random projection.  Other possible extensions of
iSVD include eigenvalue problems, linear system problems, selected
singular values within a given interval or of a given order, and
tensor decompositions. Another future direction is to explore how we
can efficiently compute the SVD if some of the columns or rows of
$\Ab$ are added (updated) or removed (downdated) after an SVD has
been obtained for a given matrix $\Ab$.

 iSVD can be accelerated using multi-level parallelism. It is
obvious the $N$ random sketches can be performed simultaneously in
parallel. The operations in each sketch and the integration process
can be parallelized as well. Efficient implementations of the
proposed algorithms on parallel computers will allow us to quickly
estimate the SVD of large-scale matrices on GPUs, parallel
computers, or distributed systems such as Spark
\cite{li2015sparkbench}.

In addition to the development of new algorithms and parallel
implementations, the tuning of parameters can affect the timing and
accuracy. Depending on the requirements (e.g., accuracy and number
of singular values), matrix structures (e.g., sparsity, size, and
distribution of the singular values), and computer architectures
(e.g., multi-core CPU or GPU cluster), we can choose between $N=1$
(rSVD) and $N>1$ (iSVD), the power exponent $q$, and the
oversampling size $p$ (and thus the dimension of the random sketches
$\ell$). Fine-tuning of Algorithm~\ref{alg:KN} or gradient-based
optimization methods may further improve the performance of iSVD.
One example is the step size used to move from the current iterate
to the next iterate.

In short, we have proposed and justified a new randomized algorithm
to compute the approximate rank-$k$ SVD of a large matrix by
integrating multiple leading subspaces based on random sketches. The
framework can be further improved and extended to benefit data
analytics, computational sciences and engineering in a broad manner.

\section*{Acknowledgments}
This work is partially supported by the Ministry of Science and
Technology, the National Center for Theoretical Sciences, and the
Taida Institute for Mathematical Sciences in Taiwan.

\bibliographystyle{plain}
\bibliography{ensemble}

\section*{Appendix}

\renewcommand{\thesubsection}{A.\arabic{subsection}}
\renewcommand{\theequation}{A.\arabic{equation}}
\setcounter{subsection}{0}

\subsection{Proof of Theorem~\ref{lemma:key}}
\label{apdnx:lemma:key}

\begin{proof}
Since $\Omegab_\i$ is a Gaussian random matrix, we have the
expectation
\begin{eqnarray*}
&& E\left(\Qb_{[i]}\Qb_{[i]}^\intercal\right)
   = E\left(\Ab\Omegab_{[i]} \left(\Omegab_{[i]}^\intercal\Ab^\intercal \Ab\Omegab_{[i]}\right)^{-1}
 \Omegab_{[i]}^\intercal\Ab^\intercal\right)\nonumber\\
&=& \Ub~ E\left(\Sigmab\Vb^\intercal\Omegab_{[i]}
  \left(\Omegab_{[i]}^\intercal\Vb\Sigmab^2\Vb^\intercal\Omegab_{[i]}\right)^{-1}
  \Omegab_{[i]}^\intercal\Vb\Sigmab\right)~ \Ub^\intercal
  =:\Ub\Lambdab\Ub^\intercal,
\end{eqnarray*}
where
\begin{equation}\label{eq:Lambda}
\Lambdab=E\left(\Sigmab\Vb^\intercal\Omegab_{[i]}
  \left(\Omegab_{[i]}^\intercal\Vb\Sigmab^2\Vb^\intercal\Omegab_{[i]}\right)^{-1}
  \Omegab_{[i]}^\intercal\Vb\Sigmab\right).
\end{equation}
Note that $\Omegab_{[i]}^\intercal\Vb\Sigmab^2
\Vb^\intercal\Omegab_{[i]}$ is non-singular with probability one.

(a)~{\it First, we show that $\Lambdab$ is a diagonal matrix.} Its
$(j,j')$th entry is given by
\[E\left(\sigma_j\sigma_{j'}\zb_j^\intercal
 \Big(\sum_{l=1}^m \sigma_l^2 \zb_l\zb_l^\intercal\Big)^{-1}
 \zb_{j'}\right),\]
where $\left[\zb_1,\dots,\zb_n \right] =
\Omegab_{[i]}^\intercal\Vb$. Note that $\zb_j
=\Omegab_{[i]}^\intercal \vb_j$, where $\vb_j$ is the $j$th column
of~$\Vb$. Let $\omegab_l$ denote the $l$th column of
$\Omegab_{[i]}$. Below we show that all off-diagonal entries of
$\Lambdab$ are zero. Without loss of generality, consider the
$(1,j)$th entry of $\Lambdab$. Let
$\Vb_{-1}:=\left[-\vb_1,\vb_2,\ldots,\vb_n\right]$ and
$\widetilde{\Omegab}_\i := \Vb \Vb_{-1}^\intercal \Omegab_\i$. Then,
$\widetilde{\Omegab}_\i^{\intercal}\Vb=\Omegab_\i^\intercal\Vb_{-1}
\Vb^\intercal \Vb = \Omegab_\i^\intercal\Vb_{-1}$. Let
$\left[\widetilde\zb_1,\dots,\widetilde\zb_n \right] :=
\widetilde{\Omegab}_\i^{\intercal}\Vb$. Then,
$\widetilde\zb_1=\widetilde{\Omegab}_\i^{\intercal}\vb_1=-\zb_1$ and
$\widetilde\zb_j=\widetilde{\Omegab}_\i^{\intercal}\vb_j=\zb_j,
~\forall j\neq 1$. Note that $\Omegab_{[i]}$ and
$\widetilde{\Omegab}_\i$ have the same distribution, as $\Omegab_\i$
have i.i.d. Gaussian entries and $(\Vb
\Vb_{-1}^\intercal)^\intercal\Vb \Vb_{-1}^\intercal=\Ib$. It implies
that $\left[\widetilde\zb_1,\dots,\widetilde\zb_n \right]$ and
$\left[\zb_1,\dots,\zb_n \right]$ follow the same distribution. That
is,
\[
\zb_1^\intercal \Big(\sum_{l=1}^m \sigma_l^2 \zb_l\zb_l^\intercal \Big)^{-1}\zb_j
\stackrel{d}{=}\widetilde\zb_1^\intercal \Big(\sum_{l=1}^m \sigma_l^2
  \widetilde\zb_l\widetilde\zb_l^\intercal\Big)^{-1}\widetilde\zb_j
= -\zb_1^\intercal
 \Big(\sum_{l=1}^m \sigma_l^2 \zb_l\zb_l^\intercal \Big)^{-1}\zb_j,
\]
where $\stackrel{d}=$ means equal in distribution. Therefore, for
the $(1,j)$th entry of $\Lambdab$, we have
\[
E\left\{\zb_1^\intercal \Big(\sum_{l=1}^m \sigma_l^2
  \zb_l\zb_l^\intercal \Big)^{-1}\zb_{j}\right\}
= -E\left\{\zb_1^\intercal \Big(\sum_{l=1}^m \sigma_l^2
  \zb_l\zb_l^\intercal \Big)^{-1}\zb_{j}\right\}= 0.
\]

(b) {\it Next, we show that all the diagonals, $E\left(\sigma_j^2
\zb_j^\intercal \Big(\sum_{l=1}^m \sigma_l^2
\zb_l\zb_l^\intercal\Big)^{-1} \zb_{j}\right)$, $j=1,\dots,n$, are
less than one.} Let $\Bb_{(-j)} := \sum_{l\neq j}^m \sigma_l^2
\zb_l\zb_l^\intercal$. As $\Omegab_\i$ consists of i.i.d. Gaussian
entries and $\ell<m$, $\Bb_{(-j)}$ is strictly positive definite
with probability one. By Sherman-Morrison-Woodbury matrix identity,
we have
\[
\Big(\sum_{l=1}^m \sigma_l^2 \zb_l\zb_l^\intercal\Big)^{-1}
=
  \Bb_{(-j)}^{-1} - \frac {\sigma_j^2  \Bb_{(-j)}^{-1} \zb_j \zb_j^\intercal
  \Bb_{(-j)}^{-1}}{1+\sigma_j^2 \zb_j^\intercal \Bb_{(-j)}^{-1} \zb_j}.
  \]
Then,
\begin{eqnarray}
&&\sigma_j^2 \zb_j^\intercal \Big(\sum_{l=1}^m \sigma_l^2 \zb_l\zb_l^\intercal\Big)^{-1} \zb_j
  = \sigma_j^2 \zb_j^\intercal \left(\Bb_{(-j)}^{-1} - \frac {\sigma_j^2  \Bb_{(-j)}^{-1} \zb_j \zb_j^\intercal
    \Bb_{(-j)}^{-1}} {1+\sigma_j^2 \zb_j^\intercal \Bb_{(-j)}^{-1} \zb_j}\right) \zb_j\nonumber\\
&=& \sigma_j^2 \left(\zb_j^\intercal \Bb_{(-j)}^{-1} \zb_j -\frac {\sigma_j^2 \zb_j^\intercal\Bb_{(-j)}^{-1}
   \zb_j \zb_j^\intercal \Bb_{(-j)}^{-1} \zb_j} {1+\sigma_j^2 \zb_j^\intercal \Bb_{(-j)}^{-1} \zb_j}\right)\nonumber\\
   &=&\frac {\sigma_j^2 \zb_j^\intercal \Bb_{(-j)}^{-1} \zb_j} {1+\sigma_j^2 \zb_j^\intercal \Bb_{(-j)}^{-1} \zb_j}
    =1-\frac 1 {1+\sigma_j^2 \zb_j^\intercal \Bb_{(-j)}^{-1} \zb_j} <1.
\label{eq:lt_one}
\end{eqnarray}
(b)~can be obtained by taking expectation of the inequality above.

(c)~{\it Finally, we want to show that $E\big\{\sigma_j^2
\zb_j^\intercal \big(\sum_{l=1}^m \sigma_l^2
\zb_l\zb_l^\intercal\big)^{-1} \zb_{j}\big\}$ is strictly decreasing
as $j$ increases.} Without loss of generality, we will only show the
comparison for $j=1,2$, i.e., $E\big\{\sigma_1^2 \zb_1^\intercal
\big(\sum_{l=1}^m \sigma_l^2 \zb_l\zb_l^\intercal\big)^{-1}
\zb_{1}\big\} > E\big\{\sigma_2^2 \zb_2^\intercal \big(\sum_{l=1}^m
\sigma_l^2 \zb_l\zb_l^\intercal\big)^{-1} \zb_{2}\big\}$. Consider
$\widetilde{\Omegab}_\i:=\Vb\Vb_{1,2}^\intercal \Omegab_\i$, where
$\Vb_{1,2}:=\left[\vb_2,\vb_1,\vb_3,\ldots,\vb_n\right]$. Let
$\left[\xb_1,\xb_2,\dots,\xb_n \right]
:=\widetilde{\Omegab}_\i^\intercal \Vb$. Then, $\xb_1=\zb_2$,
$\xb_2=\zb_1$, and $\xb_j=\zb_j$ for all $3 \leq j \leq n$. Similar
to (\ref{eq:lt_one}),
\begin{equation}\label{eq:lt_two}
\sigma_j^2 \xb_j^\intercal \Big(\sum_{l=1}^m \sigma_l^2 \xb_l\xb_l^\intercal\Big)^{-1} \xb_j
 = 1-\frac 1 {1+\sigma_j^2 \xb_j^\intercal \widetilde{\Bb}_{(-j)}^{-1} \xb_j},
\end{equation}
where $\widetilde{\Bb}_{(-j)} := \sum_{l\neq j}^m
\sigma_l^2\xb_l\xb_l^\intercal$. Again, we only need to consider the
case that $\widetilde{\Bb}_{(-j)}$ is of full rank, which holds with
probability one. Observe that
$\widetilde{\Bb}_{(-2)}=\Bb_{(-1)}+(\sigma_1^2-\sigma_2^2) \zb_2
\zb_2^\intercal$. Then,
\begin{eqnarray*}
&&\xb_2^\intercal \widetilde{\Bb}_{(-2)}^{-1} \xb_2 = \zb_1^\intercal \left(\Bb_{(-1)}
 +(\sigma_1^2-\sigma_2^2) \zb_2 \zb_2^\intercal\right)^{-1} \zb_1 \\
&=& \zb_1^\intercal \Bb_{(-1)}^{-1} \zb_1  -\frac { (\sigma_1^2-\sigma_2^2)\zb_1^\intercal\Bb_{(-1)}^{-1}
   \zb_2 \zb_2^\intercal \Bb_{(-1)}^{-1} \zb_1} {1+ \zb_{2}^\intercal \Bb_{(-1)}^{-1} \zb_{2}}
  \leq \zb_1^\intercal \Bb_{(-1)}^{-1} \zb_1.
\end{eqnarray*}
The equality holds only when $\zb_1^\intercal\Bb_{(-1)}^{-1}
\zb_2=0$, which happens with zero probability. In the following, we
will then only consider the case that $\xb_2^\intercal
\widetilde{\Bb}_{(-2)}^{-1} \xb_2 < \zb_1^\intercal \Bb_{(-1)}^{-1}
\zb_1$, which holds with probability one. Since $\sigma_1>
\sigma_2>0$, we have $\sigma_2^2 \xb_2^\intercal
\widetilde{\Bb}_{(-2)}^{-1} \xb_2<\sigma_1^2 \zb_1^\intercal
\Bb_{(-1)}^{-1} \zb_1$. Along with $\eqref{eq:lt_two}$, we have
\[
\sigma_2^2 \xb_2^\intercal \Big(\sum_{l=1}^m \sigma_l^2 \xb_l\xb_l^\intercal\Big)^{-1} \xb_2 < \sigma_1^2 \zb_1^\intercal
 \Big(\sum_{l=1}^m \sigma_l^2 \zb_l\zb_l^\intercal\Big)^{-1} \zb_1.
\]
Similarly, we have $\sigma_1^2 \xb_1^\intercal \Big(\sum_{l=1}^m
\sigma_l^2 \xb_l\xb_l^\intercal\Big)^{-1} \xb_1
  > \sigma_2^2 \zb_2^\intercal \Big(\sum_{l=1}^m \sigma_l^2 \zb_l\zb_l^\intercal\Big)^{-1} \zb_2$.
Then,
\begin{eqnarray*}
&&\sigma_1^2 \zb_1^\intercal \Big(\sum_{l=1}^m \sigma_l^2 \zb_l\zb_l^\intercal\Big)^{-1} \zb_1
  +\sigma_1^2 \xb_1^\intercal \Big(\sum_{l=1}^m \sigma_l^2 \xb_l\xb_l^\intercal\Big)^{-1} \xb_1  \\
& > & \sigma_2^2 \xb_2^\intercal \Big(\sum_{l=1}^m \sigma_l^2 \xb_l\xb_l^\intercal\Big)^{-1} \xb_2
  +\sigma_2^2 \zb_2^\intercal \Big(\sum_{l=1}^m \sigma_l^2 \zb_l\zb_l^\intercal\Big)^{-1} \zb_2.
\end{eqnarray*}
Take the expectation, and we have
\begin{eqnarray}
&&E\left(\sigma_1^2 \zb_1^\intercal \Big(\sum_{l=1}^m \sigma_l^2 \zb_l\zb_l^\intercal\Big)^{-1} \zb_1\right)
  +E\left(\sigma_1^2 \xb_1^\intercal \Big(\sum_{l=1}^m \sigma_l^2 \xb_l\xb_l^\intercal\Big)^{-1} \xb_1 \right)\nonumber \\
&>& E\left(\sigma_2^2 \xb_2^\intercal \Big(\sum_{l=1}^m \sigma_l^2 \xb_l\xb_l^\intercal\Big)^{-1} \xb_2\right)
  +E\left(\sigma_2^2 \zb_2^\intercal \Big(\sum_{l=1}^m \sigma_l^2 \zb_l\zb_l^\intercal\Big)^{-1} \zb_2.\right). \label{eq:both_increasing}
\end{eqnarray}
Since $\widetilde{\Omegab}_\i$ and $\Omegab_\i$ have the same
distribution, we have $\Omegab_\i^\intercal\Vb \stackrel{d}{=}
\widetilde{\Omegab}_\i^\intercal\Vb =\Omegab_\i^\intercal\Vb_{1,2}$,
and hence $[\zb_1,\zb_2,\zb_3,\dots,\zb_n] \stackrel{d}{=}
[\xb_1,\xb_2,\xb_3,\dots,\xb_n]$. Then,
\begin{eqnarray*}
E\left(\sigma_1^2 \zb_1^\intercal \Big(\sum_{l=1}^m \sigma_l^2 \zb_l\zb_l^\intercal\Big)^{-1} \zb_1\right)
  &=&E\left(\sigma_1^2 \xb_1^\intercal \Big(\sum_{l=1}^m \sigma_l^2 \xb_l\xb_l^\intercal\Big)^{-1} \xb_1 \right)
      \nonumber \\
E\left(\sigma_2^2 \xb_2^\intercal \Big(\sum_{l=1}^m \sigma_l^2 \xb_l\xb_l^\intercal\Big)^{-1} \xb_2\right)
  &=&E\left(\sigma_2^2 \zb_2^\intercal \Big(\sum_{l=1}^m \sigma_l^2 \zb_l\zb_l^\intercal\Big)^{-1} \zb_2.\right).
\end{eqnarray*}
Therefore,  (\ref{eq:both_increasing}) becomes
\begin{equation}\label{eq:non_incrs}
E\left\{\sigma_1^2 \zb_1^\intercal \Big(\sum_{l=1}^m \sigma_l^2
\zb_l\zb_l^\intercal\Big)^{-1} \zb_1\right\}>
E\left\{\sigma_2^2 \zb_2^\intercal \Big(\sum_{l=1}^m
\sigma_l^2 \zb_l\zb_l^\intercal\Big)^{-1} \zb_2\right\}.
\end{equation}
Similarly, we can have $E\left(\sigma_j^2 \zb_j^\intercal
\Big(\sum_{l=1}^m \sigma_l^2 \zb_l\zb_l^\intercal\Big)^{-1}
\zb_j\right)> E\left(\sigma_{j'}^2 \zb_{j'}^\intercal
\Big(\sum_{l=1}^m \sigma_l^2 \zb_l\zb_l^\intercal\Big)^{-1}
\zb_{j'}\right)$ for any pair of $(j,j')$ satisfying $j<j'$.
\end{proof}

\subsection{Proof of Theorem~\ref{lemma:q_ell}}
\label{a.proof.equiv}
\begin{proof}
Direct calculations lead to the following equalities for $\Qb\in
\St_{m,\ell}$:
\begin{eqnarray*}
  &&\left\|\Qb_\i\Qb_\i^\intercal -\Qb\Qb^\intercal\right\|_F^2
    = \trop\left\{(\Qb_\i\Qb_\i^\intercal -\Qb\Qb^\intercal)^2\right\}\\
  & =& \trop\left\{\Qb_\i\Qb_\i^\intercal\right\}
    + \trop\left\{\Qb\Qb^{\intercal}\right\}
    - 2\trop\left\{\Qb_\i\Qb_\i^\intercal\Qb\Qb^{\intercal}\right\} \\
  & =& 2\ell - 2\trop\left\{\Qb\Qb^{\intercal}\Qb_\i\Qb_\i^\intercal\right\}.
\end{eqnarray*}
Hence the summation becomes $\sum_{i=1}^N
\big\|\Qb_\i\Qb_\i^\intercal -\Qb\Qb^\intercal\big\|_F^2
  =2N\ell - 2 N \trop(\Qb\Qb^{\intercal} \overline{\Pb})$.
Similarly, one can also show that $\left\|\overline{\Pb}
-\Qb\Qb^\intercal\right\|_F^2 =\trop (\overline{\Pb}^2) + \ell - 2
\trop(\Qb\Qb^{\intercal} \overline{\Pb})$. Since $\ell$ and
$\overline{\Pb}$ are given and fixed, the two optimization problems
are equivalent.
\end{proof}

\subsection{Proof of Theorem~\ref{thm:proj_grad}}
\label{a.proof.proj_grad}

Note that the optimization in Stiefel manifold has been analyzed
in~\cite{tagare2011notes,wen2013feasible}. Here we derive the
related properties by using fundamental matrix algebras and
calculus. We hope this approach based on fundamental tools may
benefit readers who are not familiar with the advanced differential
geometry topics adopted in~\cite{tagare2011notes,wen2013feasible}.

\begin{proof}
First we find a necessary and sufficient condition for $\Xb$ being
in ${\mathcal T}_{\Qb}\St_{m,\ell}$. For all $\Xb \in {\mathcal
T}_{\Qb}\St_{m,\ell}$, find a path $\Gammab(t)$ in $\St_{m, \ell}$
with $\Gammab(0) = \Qb$ and $\Gammab'(0) = \Xb$. From
$\Gammab(t)^\intercal \Gammab(t) = \Ib$, differentiate each side by
$t$ and take $t = 0$, we have
\begin{equation}
\Xb^\intercal \Qb + \Qb^\intercal \Xb = \0,
\label{eq:tangent}
\end{equation}
which gives a necessary condition for $\Xb \in {\mathcal
T}_{\Qb}\St_{m,\ell}$. There are $\ell(\ell+1)/2$ conditions for
$\Xb$ in \eqref{eq:tangent} and the dimension of ${\mathcal
T}_{\Qb}\St_{m,\ell}$ is $m\ell - \ell(\ell+1)/2$, which means
\eqref{eq:tangent} is also a sufficient condition for $\Xb \in
{\mathcal T}_{\Qb}\St_{m,\ell}$. By taking $\vecop$ to each sides of
\eqref{eq:tangent}, we get the equality
\[
[(\Qb^\intercal  \otimes \Ib_{\ell})\Kb_{m, \ell}
+ (\Ib_\ell \otimes \Qb^\intercal )]\vecop(\Xb) = \0,
\]
where $\otimes$ denotes the Kronecker product and $\Kb_{m, \ell}$
denotes the $m \times \ell$ commutation matrix
\cite{magnus1979commutation}. Define $\Tb = \Kb_{\ell, m}(\Qb
\otimes \Ib_\ell) + (\Ib_\ell \otimes \Qb)$ and get $\Tb^\intercal
\vecop(\Xb) = \0$. This shows that the tangent space (after
vectorizing each elements) is contained in the null space of
$\Tb^\intercal$. One can compute the rank of $\Tb$ and shows that
the null space of $\Tb^\intercal$ is actually the tangent space.
Hence the projection matrix onto the tangent space is given by $(\Ib
- \Pb_{\Tb})$, where $\Pb_{\Tb} = \Tb(\Tb^\intercal
\Tb)^+\Tb^\intercal $ and $(\Tb^\intercal \Tb)^+$ denoted the
Moore-Penrose pseudo-inverse. With $\Pb_{\Tb}$, $\Db_F$ can be given
via $\vecop(\Db_F)=(\Ib - \Pb_{\Tb}) \vecop(\Gb_F)$. With some
calculation, we have $\Tb = (\Ib_\ell \otimes \Qb)(\Ib_{\ell^2} +
\Kb_{\ell,\ell})$ and thus
\begin{equation*}
\begin{split}
  \Tb^\intercal \Tb &= (\Ib_{\ell^2} + \Kb_{\ell,\ell})^\intercal (\Ib_\ell \otimes \Qb)^\intercal
  (\Ib_\ell \otimes \Qb)(\Ib_{\ell^2} + \Kb_{\ell,\ell}) \\
  & = (\Ib_{\ell^2} + \Kb_{\ell,\ell})(\Ib_{\ell^2} + \Kb_{\ell,\ell})
  = 2(\Ib_{\ell^2} + \Kb_{\ell,\ell}).
\end{split}
\end{equation*}
Then the projection matrix $\Pb_{\Tb}$ can be calculated as:
\begin{equation*}
\begin{split}
  \Pb_{\Tb} &= \Tb(\Tb^\intercal \Tb)^+\Tb^\intercal  \\
  & = (\Ib_\ell \otimes \Qb)(\Ib_{\ell^2} + \Kb_{\ell,\ell})\frac{1}{2}(\Ib_{\ell^2} + \Kb_{\ell,\ell})^+(\Ib_{\ell^2}
   + \Kb_{\ell,\ell})^\intercal (\Ib_\ell \otimes \Qb)^\intercal  \\
  & = \frac{1}{2}(\Ib_\ell \otimes \Qb)(\Ib_{\ell^2}
    + \Kb_{\ell,\ell})(\Ib_\ell \otimes \Qb^\intercal )
    \frac{1}{2}(\Ib_\ell \otimes \Qb\Qb^\intercal ) + \frac{1}{2}(\Qb^\intercal  \otimes \Qb)\Kb_{m, \ell}.
\end{split}
\end{equation*}
Hence, by $\vecop(\Db_F) = (\Ib - \Pb_{\Tb})\vecop(\Gb_F)$,
\begin{equation*}
\begin{split}
  \vecop(\Db_F)
  & = (\Ib_{\ell^2} - \frac{1}{2}(\Ib_\ell \otimes \Qb\Qb^\intercal ) - \frac{1}{2}(\Qb^\intercal  \otimes \Qb)\Kb_{m, \ell})\vecop(\Gb_F) \\
  & = \vecop(\Gb_F) - \frac{1}{2}(\Ib_\ell \otimes \Qb\Qb^\intercal)\vecop(\Gb_F)
   - \frac{1}{2}(\Qb^\intercal\otimes\Qb)\Kb_{m, \ell}\vecop(\Gb_F) \\
  & = \vecop(\Gb_F) - \frac{1}{2}\vecop(\Qb\Qb^\intercal \Gb_F) - \frac{1}{2}\vecop(\Qb\Gb_F^\intercal \Qb)
\end{split}
\end{equation*}
and $\Db_F$ can be written as
\begin{equation}
  \Db_F = \left( \Ib - \frac{1}{2}\Qb\Qb^\intercal  \right)\Gb_F - \frac{1}{2}\Qb\Gb_F^\intercal \Qb.
\label{eq:D_F}
\end{equation}
Since we have the property $\Qb^\intercal \Gb_F(\Qb) =
\Gb_F(\Qb)^\intercal \Qb$ here, we can get $\Db_F(\Qb) = (\Ib -
\Qb\Qb^\intercal)\Gb_F(\Qb)$. This completes the proof.
\end{proof}

\subsection{Proof of Lemma \ref{lemma:XX<1/4}}\label{a.proof.lemma:XX<1/4}
\begin{proof} (a)
Express $\Wb$ as $\Qb \Cb + \Qb_\bot \Bb$. Then,
$\Wb^\intercal\Wb=\Ib$ will imply $\Cb^2 + \Bb^\intercal \Bb =\Ib$.
Thus, $\Cb^4 + \Cb\Bb^\intercal \Bb\Cb =\Cb^2$. Furthermore,
$\varphi_{\Qb}(\Wb)=(\Ib-\Qb\Qb^\intercal)\Wb\Wb^\intercal\Qb =
\Qb_\bot \Bb \Cb$. Then, we have
\[
\frac{\Ib}4-\varphi_{\Qb}(\Wb)^\intercal \varphi_{\Qb}(\Wb)
=\frac{\Ib}4-\Cb \Bb^\intercal \Bb \Cb
= \frac{\Ib}4-\Cb^2 +\Cb^4 =\left(\frac{\Ib}2-\Cb^2\right)^2,
\]
which is non-negative definite. (b)~Let $\Xb_\i =
(\Ib-\Qb\Qb^\intercal)\Qb_\i\Qb_\i^\intercal \Qb$. Then, for any
vector $\vb\in\real^\ell$
\begin{eqnarray*}
&& \vb^\intercal\left(\frac{\Ib}4-\Xb^\intercal \Xb\right)\vb
 = \frac{\|\vb\|^2}4 - \Big\|\frac1N\sum_{i=1}^N \Xb_\i\vb\Big\|^2\\
&\geq& \frac{\|\vb\|^2}4-\frac 1N \sum_{i=1}^N \left\|\Xb_\i \vb\right\|^2
  = \frac 1N \sum_{i=1}^N\left(\frac 14 \|\vb\|^2- \left\|\Xb_\i \vb\right\|^2\right) \geq 0.
\end{eqnarray*}
The last inequality holds since $ \left\|\Xb_\i \vb\right\| \leq
\|\vb\|/2$ for every $i$ from~(a).
\end{proof}

\subsection{Proof of Lemma \ref{lemma:C}}\label{a.proof.lemma:C}
\begin{proof} We will show this lemma under the condition that $\Xb$ has full rank.
For $\Xb$ being rank deficient, the proof is more complicated and is
placed in Appendix~\ref{a.X_rank_deficient}. Express $\Wb=\Qb\Cb
+\Qb_\bot\Bb$. We want to find $\Bb$ and $\Cb$ satisfying
(a)~$\Wb^\intercal\Wb =\Ib_\ell$ and (b)~$\Xb=\varphi_{\Qb}(\Wb )$.
From condition~(b), it leads to $ \Xb=(\Ib-\Qb\Qb^\intercal)(\Qb
\Cb+ \Qb_\bot \Bb) (\Qb \Cb+ \Qb_\bot \Bb)^\intercal \Qb =\Qb_\bot
\Bb \Cb$. Since $\Xb$ has full rank, $\Cb$ has to have full rank and
hence is invertible. Then, $\Qb_\bot \Bb =\Xb \Cb^{-1}$. From
condition~(a), it leads to $\Ib = (\Qb \Cb+ \Qb_\bot \Bb) ^\intercal
(\Qb \Cb+ \Qb_\bot\Bb) =\Cb^2 + \Cb^{-1}\Xb^\intercal\Xb\Cb^{-1}$.
Then, $\Cb^2 = \Cb^4 + \Xb^\intercal\Xb$. With the assumption that
$\frac{\Ib}4-\Xb^\intercal\Xb$ is non-negative definite, we have
$\left\{\Cb^2-\frac{\Ib}2\right\}^2 =
\frac{\Ib}4-\Xb^\intercal\Xb,~~ \mbox{and then}~
\Cb=\left\{\frac{\Ib}2 +
\Big(\frac{\Ib}4-\Xb^\intercal\Xb\Big)^{1/2}\right\}^{1/2}$.
\end{proof}

\subsection{Matrix Square Root}\label{a.matrix_sqrt}
A matrix square root for a symmetric and non-negative definite
matrix $\Mb\in\real^{\ell\times\ell}$ is defined as follows.
\begin{equation}\label{eq:mat_sqrt}
\mbox{$\Mb^{1/2}$ is any matrix $\Tb$ that satisfies
$\Tb\Tb^\intercal=\Mb$.}
\end{equation}
Express $\Mb$ in its spectrum $\Mb=\Gb\Db\Gb^\intercal$. If we
restrict $\Tb$ to be symmetric, then
\begin{equation}\label{eq:mat_sqrt2}
\Tb = \Gb\Db^{1/2}\Gb^\intercal,~ {\rm where}~ \Db^{1/2}
=\diag\big([\pm \sqrt{d_1},\dots,\pm \sqrt{d_\ell}]\big).
\end{equation}
If $\Tb$ is further restricted to be non-negative definite, then it
is uniquely given by
\begin{equation}\label{eq:mat_sqrt3}
\Tb = \Gb\Db^{1/2}\Gb^\intercal,~ {\rm where}~ \Db^{1/2}
={\rm diag}\big([\sqrt{d_1},\dots, \sqrt{d_\ell}]\big).
\end{equation}

\subsection{Proof of Lemma~\ref{lemma:C} for Rank Deficient $\Xb$}
\label{a.X_rank_deficient}
\begin{proof}
Let $\Wb=\Qb \Cb +\Xb \Bb$. It suffices to show that $\Cb$ is
nonsingular. Then we have $\Xb \Bb =\Xb \Cb^{-1}$ again, and the
rest arguments of the proof for Lemma~\ref{lemma:C} remain the same.
If $\Cb \vb=\0 $ for some nonzero column vector $\vb$, we have
$\Xb\vb=\Xb\Bb\Cb\vb=\0$. Factorize $\Cb$ as
\[\Cb= [\Tb_1,\Tb_2]_{\ell \times \ell}\left[\begin{array}{ll}
\Db_{\ell_1\times\ell_1}& \0_{\ell_1\times(\ell-\ell_1)}\\[1ex]
\0_{(\ell-\ell_1)\times\ell_1}\ & \0_{(\ell-\ell_1)\times(\ell-\ell_1)}  \end{array}\right]\,
\left[\begin{array}{c}\Tb_1^\intercal \\[1ex]
\Tb_2^\intercal\end{array}\right]_{\ell\times\ell},
\]
where $\Db$ is diagonal and nonsingular, and
$\left[\Tb_1,\Tb_2\right]^\intercal \left[\Tb_1,\Tb_2\right]=\Ib$.
Since $\Cb \Tb_2=\0$, so is $\Xb \Tb_2$, which means that $\Xb$ can
be factorized as $\Xb = \Yb \left[ \widetilde{\Xb}_{m \times\ell_1},
\0_{m \times(\ell-\ell_1)} \right]\Tb^\intercal$ where
$\Tb=\left[\Tb_1,\Tb_2\right]$, and $\Yb^\intercal\Yb=\Ib$. Let
$\widetilde{\Bb}=\Tb^\intercal \Bb \Tb$. (b) leads to
\[
\Yb \left[\widetilde{\Xb}_{m \times\ell_1}, \0_{m \times(\ell-\ell_1)} \right]=
\Yb \left[\widetilde{\Xb}_{m \times\ell_1}, \0_{m \times(\ell-\ell_1)} \right] \widetilde{\Bb} \left[\begin{array}{ll}
\Db_{\ell_1\times\ell_1}& \0_{\ell_1\times(\ell-\ell_1)}\\[1ex]
\0_{(\ell-\ell_1)\times\ell_1}\ & \0_{(\ell-\ell_1)\times(\ell-\ell_1)}  \end{array}\right],
\]
which forces $\widetilde{\Bb}$ to be of the form $\widetilde\Bb =
\left[
\begin{array}{cc}
\widetilde\Bb_{11}& \widetilde\Bb_{12}\\[0.8ex]
\widetilde\Bb_{21}    & \widetilde\Bb_{22}
\end{array}\right]$
with $\widetilde{\Xb}_{m \times\ell_1}
\widetilde\Bb_{11}=\widetilde{\Xb}_{m \times\ell_1} \Db^{-1}$. Then
\[
\Bb^\intercal \Xb^\intercal \Xb \Bb= \Tb  \left[\begin{array}{ll}
\Db^{-1}\widetilde{\Xb}^\intercal \widetilde{\Xb}\Db^{-1}& \Db^{-1}\widetilde{\Xb}^\intercal
 \widetilde{\Xb} \widetilde{\Bb}_{12}\\[1ex]
\widetilde{\Bb}_{12}^\intercal\widetilde{\Xb}^\intercal \widetilde{\Xb} \Db^{-1} & \widetilde{\Bb}_{12}^\intercal
 \widetilde{\Xb}^\intercal \widetilde{\Xb} \widetilde{\Bb}_{12}  \end{array}\right] \Tb^\intercal.
\]
(a) leads to $\Db^{-1}\widetilde{\Xb}^\intercal \widetilde{\Xb}
\widetilde{\Bb}_{12}=\0$ and
$\widetilde{\Bb}_{12}^\intercal\widetilde{\Xb}^\intercal
\widetilde{\Xb} \widetilde{\Bb}_{12}=\Ib$, which contradicts to each
other. Therefore, $\Cb$ has to be nonsingular.
\end{proof}

\subsection{Proof of Theorem~\ref{thm:fixed_point}}\label{a.thm.fixed_point}
Since $\Omegab_\i$'s consist of all i.i.d. Gaussian entries, we have
that $\overline{\Pb}$ has distinct eigenvalues almost surely. Let
the eigenvalue decomposition of $\overline\Pb$ be denoted by $
  \overline\Pb =
  \begin{bmatrix}
    \Qb_* & \Qb_\bot
  \end{bmatrix}
  \begin{bmatrix}
    \Lambdab_1 & \0 \\
    \0 & \Lambdab_2
  \end{bmatrix}
  \begin{bmatrix}
    \Qb_* & \Qb_\bot
  \end{bmatrix}^\intercal,
$ where $\Lambdab_1 = \diag([\lambda_1, \ldots, \lambda_\ell])$ and
$\Lambdab_2 = \diag([\lambda_{\ell + 1}, \ldots, \lambda_m])$.

\begin{proof}
Suppose we start from an initial $\Qb_{\rm ini}\in{\cal
N}_\varepsilon(\Qb_*\Rb_0)$, where $\Rb_0$ is an arbitrary orthogonal
matrix and $\varepsilon$ is determined later. Denote $\Qb_0 :=\Qb_{\rm
ini}$ and $\Qb_{t+1}:= g(\Qb_t)$, where $g$ is defined
in~(\ref{fixed_point}). We can write $\Qb_0$ as $\Qb_0=\Qb_* \Rb_0
\Cb_{0,1} + \Qb_\bot \Cb_{0,2}$, where
$\Cb_{0,1}\in\real^{\ell\times\ell}$,
$\Cb_{0,2}\in\real^{(m-\ell)\times\ell}$ and $\Cb_{0,1}^\intercal
\Cb_{0,1}+ \Cb_{0,2}^\intercal \Cb_{0,2} =\Ib_\ell$. We first show
that $\Rb_{0,*}^\intercal \Rb_0 \Cb_{0,1}$ is symmetric, where
\[
\Rb_{0,*}:=\argmin_{\Rb\in\St_{\ell,\ell}} \|\Qb_0-\Qb_*\Rb\|_F
=\argmin_{\Rb\in\St_{\ell,\ell}} \|\Qb_* \Rb_0 \Cb_{0,1}-\Qb_*\Rb\|_F
=\argmin_{\Rb\in\St_{\ell,\ell}} \|\Cb_{0,1}- \Rb_0^\intercal \Rb\|_F.
\]
Denote the spectrum of $\Cb_{0,1}$ as $\Cb_{0,1}=\Lb \Sb
\Hb^\intercal$, where the diagonal entries of $\Sb$ are less than or
equal to one. Therefore, the best approximation is given by
$\Rb_0^\intercal\Rb_{0,*}=\Lb \Hb^\intercal$. Then,
$\Rb^\intercal_{0,*} \Rb_0\Cb_{0,1}= \Hb \Lb^\intercal \Lb \Sb
\Hb^\intercal=\Hb \Sb \Hb^\intercal$. This completes the proof of
symmetry. With this symmetry property, we can write $\Qb_0$ as
$\Qb_0=\Qb_* \Rb_{0, *} \Cb_{1} + \Qb_\bot \Cb_{2}$ for some
symmetric $\Cb_{1}\in\real^{\ell\times\ell}$ and some
$\Cb_{2}\in\real^{(m-\ell)\times\ell}$.

Next, we (i) extend the definition of $g$ to an open covering ${\cal
O}_{\varepsilon_0}$ of the compact Stiefel manifold $\St_{m\times\ell}$
for some small $\varepsilon_0>0$. (ii) Compute its derivative (the
usual derivative in the Euclidean space) $\Db(\Qb_{*}\Rb_{0,*})
:={\partial \vec(g)}/{\partial\vec(\Mb)^\intercal}
\big|_{\Mb=\Qb_{*}\Rb_{0,*}}$, where $\Mb\in{\cal O}_{\varepsilon_0}$.
(iii)~Show that $\Db$'s spectral norm, when restricted to the
subspace ${\cal V}_0$ given below~(\ref{subspaceV}), is strictly
less than one. (iv) $\norm{\vec(\Qb_1) - \vec(\Qb_{*}\Rb_{0,*})}_2
\leq \sqrt{\alpha} \norm{\vec(\Qb_0) - \vec(\Qb_{*}\Rb_{0,*})}_2$
for some $0\le\alpha<1$. (v)~$\norm{\vec(\Qb_1) -
\vec(\Qb_{*}\Rb_{1,*})}_2 \leq \sqrt{\alpha} \norm{\vec(\Qb_0) -
\vec(\Qb_{*}\Rb_{0,*})}_2$, where $\Rb_{1,*}$ is defined
below~(\ref{eq:R1star}). (vi)~Iteratively obtain
$\norm{\vec(\Qb_{t+1}) - \vec(\Qb_{*}\Rb_{t+1,*})}_2 \leq
\alpha^{(t+1)/2} \norm{\vec(\Qb_0) - \vec(\Qb_{*}\Rb_{0,*})}_2$.
(vii)~Finally, establish the convergence of $\Qb_t$.

(i) For $\varepsilon_0$ sufficiently small, matrices $\Cb(\Qb)$ and
$\Xb(\Qb)$ as functions of $\Qb$ can be extended to an open covering
${\cal O}_{\varepsilon_0} :=\{\Mb\in\real^{m\times\ell}:
\inf_{\Qb\in\St_{m,\ell}} \|\Mb-\Qb\|_F<\varepsilon_0\}$ of the compact
Stiefel manifold $\St_{m\times\ell}$. Then, $g$ can be extended as
well. From now on, we consider $g$ as a function defined on this
open covering ${\cal O}_{\varepsilon_0}$ and we can take derivative of
$g$ in the usual Euclidean sense.

(ii) Recall $\Xb(\Qb) =(\Ib -\Qb\Qb^\intercal)\overline{\Pb}\Qb$,
$\Cb= \left\{\frac{\Ib}2 + \left(\frac{\Ib}4
-\Xb^\intercal\Xb\right)^{1/2}\right\}^{1/2}$, and $\Cb$ satisfies
the equation $\Cb^4-\Cb^2 +\Xb^\intercal\Xb=\0$. By taking
derivative for both sides of the last equation, we have
$\frac{\partial\vec(\Cb^4-\Cb^2 +\Xb^\intercal\Xb)}
{\partial\vec(\Qb)^\intercal} \big|_{\Qb_*\Rb_{0,*}}
 =\0$. The left side of the equation
can be calculated as follows.
\begin{eqnarray*}
&& \left\{\Cb^3\otimes\Ib + \Cb^2\otimes\Cb +\Cb\otimes\Cb^2 +
  \Ib\otimes\Cb^3 -\Cb\otimes\Ib -\Ib\otimes\Cb\right\}
  \frac{\partial\vec(\Cb)}{\partial\vec(\Qb)^\intercal}\\
&& +\left\{\Ib\otimes\Xb^\intercal+(\Xb^\intercal\otimes\Ib)\Kb_{m,\ell}\right\}
   \frac{\partial\vec(\Xb)}{\partial\vec(\Qb)^\intercal}~ \bigg|_{\Qb_*\Rb_{0,*}}
 =  \frac{2\partial\vec(\Cb)}{\partial\vec(\Qb)^\intercal} \bigg|_{\Qb_*\Rb_{0,*}},
\end{eqnarray*}
where the equality holds by the facts $\Xb(\Qb_*\Rb_{0,*}) = \0$ and
$\Cb(\Qb_*\Rb_{0,*}) = \Ib$. Hence, we have
$\frac{\partial\vec(\Cb)}{\partial\vec(\Qb)^\intercal}
\big|_{\Qb_*\Rb_{0,*}}=\0$. Applying similar techniques to both
sides of the equation $g(\Qb)\Cb = \Qb\Cb^2 + \Xb$, we get
\[
  \frac{\partial \vec(g)}{\partial\vec(\Qb)^\intercal} \bigg|_{\Qb_*\Rb_{0,*}}
  = \Ib_{m\ell} + \frac{\partial \Xb}{\partial\vec(\Qb)^\intercal} \bigg|_{\Qb_*\Rb_{0,*}}.
\]
Some further calculation goes as follows.
\begin{eqnarray*}
\Db(\Qb_*\Rb_{0, *})  &=& \Ib_{m\ell} + \frac{\partial \vec(\Xb)}{\partial \vec(\Qb)^\intercal} \bigg|_{\Qb_*\Rb_{0,*}}
   = \Ib_{m\ell} + \frac{\partial \vec((\Ib-\Qb\Qb^\intercal)\overline\Pb\Qb)}
   {\partial \vec(\Qb)^\intercal} \bigg|_{\Qb_*\Rb_{0,*}}\\
&=& \Ib_{m\ell}+\Ib_\ell\otimes\overline{\Pb} -
    \Qb^\intercal\overline{\Pb}\Qb\otimes\Ib_m -
     (\Qb^\intercal\overline{\Pb}\otimes\Qb)\Kb_{m,\ell}
     - \Ib_\ell\otimes \Qb\Qb^\intercal\overline{\Pb}~ \bigg|_{\Qb_*\Rb_{0,*}}\\
&=& \Ib_{m\ell}+ \Ib_\ell\otimes\overline{\Pb} -
   \Rb_{0,*}^\intercal\Lambdab_1\Rb_{0,*}\otimes\Ib_m -
     (\Rb_{0,*}^\intercal\Lambdab_1\Qb_*^\intercal\otimes\Qb_*\Rb_{0,*})\Kb_{m,\ell}
     - \Ib_\ell\otimes(\Qb_*\Lambdab_1\Qb_{*}^\intercal).
\end{eqnarray*}
Let
\begin{equation}\label{subspaceV}
{\cal V}_0:=\{\Vb: \Vb = \Qb_*\Rb_{0,*} \Sb_1 + \Qb_\bot \Sb_2~~
\mbox{with symmetric $\Sb_1$}\}.
\end{equation}
For any $m \times \ell$ matrix $\Vb\in{\cal V}_0$, we have
\begin{equation*}
\begin{split}
  \Db(\Qb_*\Rb_{0, *})\,  \vec(\Vb) &= \vec\left(\Vb + \overline{\Pb}{\Vb} - \Vb \Rb_{0,*}^\intercal \Lambdab_1\Rb_{0,*}
    - \Qb_{*}\Rb_{0,*}\Vb^\intercal \Qb_{*}\Lambdab_1\Rb_{0,*}
    -\Qb_{*}\Lambdab_1\Qb_{*}^\intercal\Vb \right) \\
  &= \vec\left(\Vb + \Qb_\bot\Lambdab_2\Qb_\bot^\intercal\Vb
    - \Vb\Rb_{0,*}^\intercal\Lambdab_1 \Rb_{0,*}
    - \Qb_{*}\Rb_{0,*}\Vb^\intercal \Qb_{*}\Lambdab_1\Rb_{0,*} \right) \\
  &= \vec\left(\Qb_*\Rb_{0,*}(\Sb_1 - \Sb_1\Rb_{0,*}^\intercal\Lambdab_1 \Rb_{0,*}
    - \Sb_1^\intercal\Rb_{0,*}^\intercal\Lambdab_1 \Rb_{0,*})\right)\\
  &~~~~ +\vec\left( \Qb_\bot(\Sb_2 + \Lambdab_2\Sb_2 - \Sb_2\Rb_{0,*}^\intercal\Lambdab_1 \Rb_{0,*})\right).
\end{split}
\end{equation*}

(iii) Then,
\begin{equation*}
\begin{split}
  \norm{\Db(\Qb_*\Rb_{0, *})\, \vec(\Vb)}_2^2 &=
  \norm{\Qb_{*}\Rb_{0,*}(\Sb_1 - 2 \Sb_1 \Rb_{0,*}^\intercal \Lambdab_1 \Rb_{0,*} ) + \Qb_\bot(\Sb_2 + \Lambdab_2\Sb_2 - \Sb_2 \Rb_{0,*}^\intercal \Lambdab_1 \Rb_{0,*})}_F^2 \\
  &=    \norm{\Sb_1 - 2 \Sb_1 \Rb_{0,*}^\intercal \Lambdab_1  \Rb_{0,*} }_F^2 + \norm{\Sb_2 + \Lambdab_2\Sb_2 - \Sb_2\Rb_{0,*}^\intercal \Lambdab_1\Rb_{0,*}}_F^2 \\
  &\leq  (1-2\lambda_\ell)^ 2 \norm{\Sb_1}_F^2 + (1-\lambda_\ell + \lambda_{\ell+1})^2\norm{\Sb_2}_F^2 \\
  &\leq \beta \norm{\Sb_1}_F^2 + \beta \norm{\Sb_2}_F^2 =\beta \norm{\vec(\Vb)}_2^2,
\end{split}
\end{equation*}
where $\beta=\max\{(1-2\lambda_\ell)^2, (1-\lambda_\ell +
\lambda_{\ell+1})^2\}<1$. Let $\alpha = \frac{1+\beta}{2}$. By
continuity of $\Db$, there exists an $\varepsilon\in (0,\varepsilon_0)$
such that $\norm{\Db(\Mb) \vec(\Vb)}_2^2 \leq
\alpha\norm{\vec(\Vb)}_2^2$ for any $\Mb\in{\cal
B}_\varepsilon(\Qb_*\Rb_{0, *})$, where ${\cal B}_\varepsilon(\Qb_*\Rb_{0,
*})$ is an $\varepsilon$-open ball in Euclidean space, and for any
$\Vb \in{\cal V}_0$. The selection of $\varepsilon$ can be made
independent of $\Rb_{0,*}$ due to the fact that the underlying
matrix Stiefel manifold is compact and it can be covered by finitely
many $\varepsilon$-ball for any given $\varepsilon$.

(iv)~Consider a path connecting $\Qb_*\Rb_{0,*}$ and $\Qb_0$:
\[
\Mb(\tau)=\tau(\Qb_*\Rb_{0,*}\Cb_1+\Qb_\bot \Cb_{2})
  + (1-\tau) \Qb_* \Rb_{0,*},
\]
which is the line segment between $\Qb_0$ and $\Qb_*\Rb_{0,*}$ with
$\Mb(1)=\Qb_0$ and $\Mb(0)=\Qb_*\Rb_{0,*}$. We have $\Mb'(\tau) =
\Qb_*\Rb_{0,*} \Cb_{1} - \Qb_* \Rb_{0,*} + \Qb_\bot \Cb_{2}\in{\cal
V}_0$. By Mean Value Theorem on this curve,
\begin{equation*}
\begin{split}
  &\norm{\vec(\Qb_{1}) - \vec(\Qb_{*}\Rb_{0,*})}_2 =
    \norm{\vec(g(\Mb(1))) - \vec(g(\Mb(0)))}_2\\
  \leq~& \sup_{\tau \in [0, 1]} \norm{\Db_{\Mb(\tau)} \vec(\Mb'(\tau))}_2
    \leq \sqrt{\alpha} \norm{\vec(\Qb_0) - \vec(\Qb_{*}\Rb_{0,*})}_2.
\end{split}
\end{equation*}

(v)~We can define $\Rb_{1,*}$ in a similar way as how $\Rb_{0,*}$ is
defined:
\begin{equation}\label{eq:R1star}
\Rb_{1,*}:=\argmin_{\Rb\in\St_{\ell,\ell}} \|\Qb_1-\Qb_*\Rb\|_F.
\end{equation}
Note that $\Qb_1$ is closer to $\Qb_{*}\Rb_{1,*}$ than to
$\Qb_{*}\Rb_{0,*}$. Then,
\[
\norm{\vec(\Qb_{1}) - \vec(\Qb_{*}\Rb_{1,*})}_2
\leq \norm{\vec(\Qb_{1}) - \vec(\Qb_{*}\Rb_{0,*})}_2
\leq \sqrt{\alpha} \norm{\vec(\Qb_0) - \vec(\Qb_{*}\Rb_{0,*})}_2.
\]

(vi)~Furthermore, we can write $\Qb_1$ as $\Qb_1=\Qb_* \Rb_0
\Cb_{1,1} + \Qb_\bot \Cb_{1,2}$, where
$\Cb_{1,1}\in\real^{\ell\times\ell}$,
$\Cb_{1,2}\in\real^{(m-\ell)\times\ell}$ and $\Cb_{1,1}^\intercal
\Cb_{1,1}+ \Cb_{1,2}^\intercal \Cb_{1,2} =\Ib_\ell$. Following the
same arguments for $\Rb_{0,*}$, we have $\Rb_{1,*}^\intercal \Rb_0
\Cb_{1,1}$ is symmetric. That is, $\Qb_1$ can be expressed as $\Qb_1
= \Qb_{*}\Rb_{1,*} \left(\Rb_{1,*}^\intercal
 \Rb_0 \Cb_{1,1}\right)+\Qb_\bot \Cb_{1,2}$.
By similar arguments as in (iii)-(iv), we now work on
$\Db(\Qb_*\Rb_{1,*})\vec(\Vb)$, where $\Vb\in{\cal V}_1:= \{\Vb: \Vb
= \Qb_*\Rb_{1,*} \Sb_1 + \Qb_\bot \Sb_2~~ \mbox{with symmetric
$\Sb_1$}\}$, and similar derivations lead to
\[
\norm{\vec(\Qb_{2}) - \vec(\Qb_{*}\Rb_{2,*})}_2  \leq \sqrt{\alpha}
 \norm{\vec(\Qb_1) - \vec(\Qb_{*}\Rb_{1,*})}_2.
\]
Iteratively, we have
\[\begin{split}
&~\norm{\vec(\Qb_{t+1}) - \vec(\Qb_{*}\Rb_{t+1,*})}_2  \leq \sqrt{\alpha}
 \norm{\vec(\Qb_t) - \vec(\Qb_{*}\Rb_{t,*})}_2 \\
\leq&~ \alpha^{\frac {t+1} 2} \norm{\vec(\Qb_0) - \vec(\Qb_{*}\Rb_{0,*})}_2.
\end{split}
\]

(vii)~Therefore, $ \norm{\vec(\Qb_t) - \vec(\Qb_{*}\Rb_{t,*})}_2$
converges to 0. This implies that $ \Qb_t \Qb_t ^\intercal$
converges to $\Qb_{*}\Rb_{t,*} \Rb_{t,*}^\intercal \Qb_{*}^\intercal
= \Qb_* \Qb_* ^\intercal$, which is independent of $t$. It further
implies that $\Xb(\Qb_t)$ converges to $\0$, and then $\Cb(\Qb_t)$
to $\Ib_\ell$. Finally, this leads to the convergence of $\Qb_t$.
\end{proof}

\subsection{Proof of Lemma~\ref{lemma:approx_ell} and Theorem~\ref{thm:cnstcy}}\label{proof.consistency}

Proof of Lemma~\ref{lemma:approx_ell} is given below.
For $\Qb\in\St_{m,\ell}$, we have
$\| \Ub\Lambdab\Ub^\intercal - \Qb \Qb^\intercal\|_F^2 = \|
\Lambdab-\Ub^\intercal \Qb\Qb^\intercal\Ub \|_F^2$. Because
$\Ub^\intercal \Qb \in\St_{m,\ell}$, $\Ub^\intercal
\Qb\Qb^\intercal\Ub$ is a rank-$\ell$ projection matrix. The best
rank-$\ell$ projection matrix to approximate $\Lambdab$ is
$\Ub^\intercal \Qb_\opt \Qb_\opt^\intercal\Ub  = \left[
\begin{array}{cc}
\Ib_{\ell}&\0\\
\0& \0
\end{array}
\right]. $ This fact suggests that $\Qb_\opt \Qb_\opt^\intercal =
\Ub_{\ell}\Ub_{\ell}^\intercal$.

Proof of Theorem~\ref{thm:cnstcy} is given below.
(a) From Theorem~\ref{thm:lln}, we know that $\overline{\Pb}
=\frac1N\sum_{i=1}^N \Qb_{[i]}\Qb_{[i]}^\intercal$ converges to
$\Ub\Lambdab\Ub^\intercal$ with probability one. By
Theorem~\ref{lemma:q_ell} and Lemma~\ref{lemma:approx_ell}, we have
that $\overline{\Qb}\,\overline{\Qb}^\intercal$ converges to
$\Ub_{\ell} \Ub_{\ell}^\intercal$ with probability one. (b) Because
$\overline{\Qb}\,\overline{\Qb}^\intercal \to \Ub_{\ell}
\Ub_{\ell}^\intercal$ with probability one, we have $
\overline{\Qb}\,\overline{\Qb}^\intercal \Ab \to \Ub_{\ell}
\Ub_{\ell}^\intercal \Ab = \Ub_{\ell}\, \Sigmab_{\ell}\,
\Vb_{\ell}^\intercal $ with probability one. Note that
$\overline{\Qb}^\intercal\Ab = \widehat{\Wb}_{\ell}
\widehat{\Sigmab}_{\ell} \widehat{\Vb}_{\ell}^\intercal$, and
$\overline{\Qb}\,\overline{\Qb}^\intercal\Ab =
\overline{\Qb}\widehat{\Wb}_{\ell} \widehat{\Sigmab}_{\ell}
\widehat{\Vb}_{\ell}^\intercal =\widehat{\Ub}_{\ell}
\widehat{\Sigmab}_{\ell} \widehat{\Vb}_{\ell}^\intercal$.
Specifically, we have $\widehat{\Ub}_{\ell} \widehat{\Sigmab}_{\ell}
\widehat{\Vb}_{\ell}^\intercal \to \Ub_{\ell}\, \Sigmab_{\ell}\,
\Vb_{\ell}^\intercal$ with probability one. By the continuity of
left singular vectors as functions of the matrix $\Ab$, which has
distinct leading $\ell$ singular values, we have that
$\widehat{\ub}_j$ converges with probability one to $\ub_j$ up to a
sign change. Specifically, $\left|\widehat{\ub}_j^\intercal
\ub_j\right|\to 1$.

\subsection{Proof of Theorem~\ref{thm:clt2}}\label{proof.clt}

Denote $\widehat \ub_j$ as $\widehat \ub_j = h_j\big(\frac1N
\sum_{i=1}^N \Qb_{[i]}\Qb_{[i]}^\intercal\big)$. Then, by SLLN
\[\ub_j =\lim_{N\to\infty}h_j\big(\frac1N \sum_{i=1}^N
\Qb_{[i]}\Qb_{[i]}^\intercal\big) = h_j
\big(\lim_{N\to\infty}\frac1N \sum_{i=1}^N
\Qb_{[i]}\Qb_{[i]}^\intercal\big) = h_j(\Ub\Lambdab\Ub^\intercal).\]
To apply the delta-method, we need to compute the derivative
$\Deltab_j$. Let $\Mb := \Ub\Lambdab\Ub^\intercal$. Since $\ub_j$ is
the $j$th eigenvector, we have $\Mb \ub_j =\lambda_j \ub_j$, where
$\ub_j^\intercal\ub_j=1$. Let $\dot\Mb$ denote a small perturbation
to $\Mb$, and $\dot\ub_j$ and $\dot\lambda_j$ be corresponding
perturbations. Consider small perturbations to both sides of the
equation above. Then,
\[\dot\Mb \ub_j + \Mb \dot\ub_j
  = \dot\lambda_j \ub_j + \lambda_j \dot\ub_j,~~
{\rm where}~\ub_j^\intercal\dot\ub_j =0.\]
Rearrange the equation above, and we have
\begin{equation}
\left(\lambda_j\Ib_m -\Mb\right)\dot\ub_j
  = \dot\Mb\ub_j -\dot\lambda_j \ub_j.
  \label{eq:perturb}
\end{equation}
Let $\left(\lambda_j\Ib_m -\Mb\right)^+$ be the Moore-Penrose pseudo
inverse. Multiply it to both sides of Equation~(\ref{eq:perturb}), we have
$ \dot\ub_j = \left(\lambda_j\Ib_m -\Mb\right)^+
    \left(\dot\Mb\ub_j -\dot\lambda_j \ub_j\right)
    = \left(\lambda_j\Ib_m -\Mb\right)^+ \dot\Mb\ub_j.
$
Then, $ \dot\ub_j= \big[\ub_j^\intercal\otimes \left(\lambda_j\Ib_m
-\Mb\right)^+\big] \vec(\dot\Mb)$. Therefore,
\begin{eqnarray}\label{Delta}
\Deltab_j &=& \frac{\partial \ub_j}{\partial\vec
(\Mb)^\intercal} = \ub_j^\intercal\otimes
\left(\lambda_j\Ib_m -\Mb\right)^+
 = \ub_j^\intercal\otimes
\left(\lambda_j\Ib_m -\Ub\Lambdab\Ub^\intercal\right)^+ .
\end{eqnarray}
The asymptotic normality can be obtained by a straightforward
application of the delta-method.

\end{document}